\documentclass[12pt]{article}
\usepackage{amsmath,amssymb,amsfonts,amsthm}
\usepackage{latexsym}
\usepackage[english]{babel}
\usepackage{showidx}
\usepackage{graphicx}
\usepackage{graphics}
\usepackage{babel,amsmath,amssymb,amsbsy,amsfonts,latexsym,amsthm,mathrsfs}

\def\<{\langle}
\def\>{\rangle}

\def\I{\mbox{\large \bf 1}}
\def\P{{\mathbb P}}
\def\R{{\mathbb R}}
\def\E{{\mathbb E}}
\def\D{{\mathbb D}}
\def\N{{\mathbb N}}
\def\leb{\mathrm{Leb}}
\def\si{{\gamma}}

\def\cl#1{{\mathcal #1}}

\numberwithin{equation}{section}
\newtheorem{theorem}{Theorem}[section]

\newtheorem{example}[theorem]{Example}

\newtheorem{lemma}[theorem]{Lemma}

\newtheorem{proposition}[theorem]{Proposition}

\begin{document}

\title{ \textbf{Positivity and lower bounds\\ for the density of Wiener functionals}}
\author{ \textsc{Vlad Bally}\thanks{%
Laboratoire d'Analyse et de Math\'ematiques Appliqu\'ees, UMR 8050,
Universit\'e Paris-Est Marne-la-Vall\'ee, 5 Bld Descartes, Champs-sur-Marne,
77454 Marne-la-Vall\'ee Cedex 2, France. Email: \texttt{bally@univ-mlv.fr} }%
\smallskip \\
\textsc{Lucia Caramellino}\thanks{%
Dipartimento di Matematica, Universit\`a di Roma - Tor Vergata, Via della
Ricerca Scientifica 1, I-00133 Roma, Italy. Email: \texttt{%
caramell@mat.uniroma2.it}}\smallskip\\
}
\maketitle

\parindent 0pt

{\textbf{Abstract.}}
We consider a functional on the Wiener space which is
smooth and not degenerated in Malliavin sense and we give a criterion for
the strict positivity of the density, that we can use  to state lower bounds as well.
The results are based on the representation of the density in terms of the
Riesz transform introduced in Malliavin and Thalmaier \cite{bib:[M.T]} and on the
estimates of the Riesz transform given in Bally and Caramellino \cite{bib:[B.C]}.

\medskip

{\textbf{Keywords}}: Riesz transform,
Malliavin calculus, strict positivity and lower bounds for the density.

\medskip

{\textbf{2000 MSC}}: 60H07, 60H30.


\section{Introduction}
The aim of this paper is to study the strict positivity and
lower bounds for the density of a functional on the
Wiener space. Although the two problems are related each other,
the hypotheses under which the results may be obtained are
different. Just to make clear what we expect to be these
hypotheses, consider the example of a $d$ dimensional diffusion
process $X_{t}$ solution to $dX_{t}=\sum_{j=1}^{m}\sigma
_{j}(X_{t})\circ dW_{t}^{j}+b(X_{t})dt$ where $\circ dW_{t}^{j}$
denotes the Stratonovich integral. The skeleton associated to this
diffusion process is the solution $x_{t}(\phi )$ to the equation
$dx_{t}(\phi )=\sum_{j=1}^{m}\sigma _{j}(x_{t}(\phi ))\phi
_{t}^{j}dt+b(x_{t}(\phi ))dt$, for a square integrable $\phi$.
The celebrated support theorem of
Stroock and Varadhan guarantees that the support of the law of
$X_{t}$ is the closure of the set of points $x$ which are
attainable by a skeleton, that is $x=x_{t}(\phi )$ for some
control $\phi\in L^2([0,T]) .$ Suppose now that the law of $X_{t}$ has a
continuous density $p_{X_t}$ with respect to the Lebesgue
measure. Then in order to get a criterion for $p_{X_t}(x)>0$,
we prove that this holds if $x$ is attainable, that is $x=x_{t}(\phi
)$ for some $\phi$, and a suitable non degeneracy assumption holds in $x.$ The second
problem is to give a lower bound for $p_{X_t}(x)$ and this can be achieved
if a non degeneracy condition holds all along
the curve $x(\phi )$ which arrives in $x$ at time $t$. Roughly speaking,
the idea is the following. 
If one has a non degeneracy
condition all along the skeleton curve arriving in $x$ at time $t$, one may give a lower bound for the probability to remain in the tube up to $t-\delta $ for a
small $\delta >0$ and then one employs an argument based on
Malliavin calculus in order to focus on the point $x$ -
essentially this means that one is able to give a
precise estimate of the behavior of the diffusion in short time (between $%
t-\delta $ and $t$). This allows one to obtain a lower bound for $p_{X_t}(x).$ If
one is not interested in lower bounds but only in the strict positivity property, the argument is the same but one does not need to estimate the probability to remain in
the tube: using the support theorem one knows that this probability is
strictly positive (but this is just qualitative, so one has no lower bound
for it) and then one focuses on the point $x$ using again the same argument
concerning the behavior of the diffusion in short time. So one needs the non
degeneracy condition in $x$ only.

\smallskip

The two problems mentioned above have been intensively studied in the
literature. Let us begin with the strict positivity. At the best of our
knowledge the first probabilistic approach to this problem is due to
Ben Arous and L\'eandre \cite{bib:[BA.L]}, who used Malliavin calculus in order to give
necessary and sufficient conditions allowing one to have $p_{X_t}(x)>0$ for a
diffusion process (as above). They proved that if H\"{o}rmander's condition
holds then $p_{X_t}(x)>0$ if and only if $x$ is attainable by a skeleton $x_{t}(\phi )$ such that $\psi\mapsto x_{t}(\psi)$ is a submersion in $%
\phi .$ The argument they used is based on the inverse function theorem and
on a Girsanov transformation. All the papers which followed developed in some way their techniques. First, Aida, Kusuoka and Stroock \cite{bib:[A.K.S]} gave a
generalization of this criterion in an abstract framework which still
permits to exhibit a notion of skeleton. Then Hirsch and Song \cite{bib:[H.S]} studied
a variant of such criterion for a general functional on the Wiener space
using capacities and finally L\'eandre \cite{bib:[L2005]} obtained similar results for
diffusion processes on manifolds. Notice that once one has a criterion of
the above type there is still a non trivial problem to be solved: one has to
exhibit the skeleton which verifies the submersion property. So, number of
authors dealt with concrete examples in which they are able to use in a more
or less direct way the argument of Ben Arous and L\'eandre: Bally and Pardoux
\cite{bib:[B.P]} dealt with parabolic stochastic heat equations, Millet and
Sanz-Sol\'e \cite{bib:[M.SS2]} worked with hyperbolic stochastic partial differential
equations, Fournier \cite{bib:[F]} considered  jump type equations, Dalang and
E. Nualart \cite{bib:[D.N]} used such positivity results for building a potential
theory for SPDE's and E. Nualart \cite{bib:EN} has recently proved results in this direction again for solutions to SPDE's.

Concerning lower bounds for the density,
a first result was found by Kusuoka and Stroock \cite{bib:[KS3]} for diffusion
processes that verify a strong uniform H\"{o}rmander condition. Afterwards
Kohatsu-Higa \cite{bib:[K-H]} obtained lower bounds for general functionals on the
Wiener space under a uniform ellipticity condition
and Bally \cite{bib:bally} proved results under local ellipticity conditions.
Recently, Gaussian type lower and upper bounds are studied in E. Nualart and Quer-Sardanyons \cite{bib:EN2} for
the nonlinear stochastic heat equation.

\smallskip

The present  paper gives a contribution in this framework: we study
the strict positivity and lower bounds for the density of a general functional on the Wiener space starting from a result
(Proposition \ref{exp-est}) which gives the behavior of a
small perturbation of a
Gaussian random variable - it corresponds to the study of a diffusion
process in short time (between $%
t-\delta $ and $t$). This is a consequence of an abstract result (Theorem \ref{th-dist}) in which the distance between the local density functions of two random variables (doesn't matter if one of them is Gaussian) is studied. It is worth to stress that
Theorem \ref{th-dist}
is of interest in itself and can be linked to the implicit function theorem in order to get further estimates which can be used to handle the same problem under H\"ormander type conditions (see  \cite{bib:tubes}).

So, our main result (see Theorem \ref{th-pos}) gives sufficient conditions in order to obtain the following lower bound for the law of $F$ around a point $y\in \R^d$: there exists $\eta>0$ and $c(y)>0$ such that
$$
\P(F\in A)\geq c(y)\leb_d(A)\quad \mbox{for every Borel set } A\subset B_\eta(y),
$$
$\leb_d$ denoting the Lebesgue measure on $\R^d$.
In particular, if the law of $F$ is absolutely continuous on $B_\eta(y)$ then the density $p_F$ satisfies $p_F(x)\geq c(y)>0$ for every $x\in B_\eta(y)$. Essentially, our conditions are that $y$ belongs to the support of the law of $F$ and an ellipticity-type condition holds around $y$.

\smallskip

In our examples, we first deal with an It\^o process $X_t$ defined as a component of a diffusion
process, that is
\begin{align*}
X_t=
&
x_0+\sum_{j=1}^m\int_0^t\sigma_j(X_t,Y_t)dW^j_t+\int_0^tb(X_t,Y_t)dt\\
Y_t=
&
y_0+\sum_{j=1}^m\int_0^t\alpha_j(X_t,Y_t)dW^j_t+\int_0^t\beta(X_t,Y_t)dt.
\end{align*}
Notice that for diffusion processes, we get an example which
is essentially the same treated in Ben Arous and L\'eandre \cite{bib:[BA.L]} and in Aida, Kusouka
and Stroock \cite{bib:[A.K.S]}.
Let $(x(\phi),y(\phi))$ denote the skeleton associated to the diffusion pair $(X,Y)$
and let $x=x_t(\phi)$ for some suitable control $\phi$.
Then, whenever a continuous local density $p_{X_t}$ of $X_t$ exists in $x$,
we prove that if $\sigma\sigma^*(x,y_t(\phi))>0$ then
$p_{X_t}(x)>0$.
And moreover, if $\inf_{s\leq t}\inf_y\sigma\sigma^* \big(x_s(\phi),y\big)\geq \lambda_*>0$
and $x_s(\phi)$ belongs to a suitable class of paths (see Theorem \ref{th-Ito} for details), then a lower bound for $p_{X_t}(x)$ can be written in terms
of the lower estimates for the probability that It\^o processes remain near a path proved in
Bally, Fern\'andez and Meda in \cite{bfm}.

As a second example, in Section \ref{asian} we treat the two dimensional diffusion process
\[
dX_{t}^{1}=\sigma_1 (X_{t})dW_{t}+b_{1}(X_{t})dt,\quad
dX_{t}^{2}=b_{2}(X_{t})dt
\]
which is degenerated in any point $x\in \R^{2}.$ We assume that $x$ is
attainable by a skeleton $x_{t}(\phi )$ and that $|\sigma_1 (x)|>0$ and $%
|\partial _{1}b_{2}(x)|>0$ - which amounts to say that the weak H\"{o}rmander condition
holds in the point $x$. We
prove that under this hypothesis one has $p_{X_t}(x)>0.$ For this example Bally and
Kohatsu-Higa \cite{bib:[B.KH]} have already given a lower bound for the density  under the
stronger hypothesis that $\inf_{s\leq t}|\sigma (x_{s}(\phi ))|>0$ and $%
\inf_{s\leq t}|\partial _{1}b_{2}(x_{s}(\phi ))|>0$. So the same non degeneracy
condition holds but along the whole  curve $x_{s}(\phi ),0\leq s\leq t.$ Notice that
we use a skeleton $x_{s}(\phi )$ which arrives in $x$ but we do not ask for the
immersion property (according to the result of Ben Arous and L\'eandre it follows
that a skeleton which verifies the immersion property exists also, but we do
not know how to produce it directly and we do not need it). And it seems clear
to us that our criterion may be used for SPDE's as well and would
simplify the proofs given in the already mentioned papers.

\smallskip

The paper is organized as follows. In Section \ref{sect-resume}
we first state localized representation formulas for the density by means of the Riesz
transform (see Section \ref{sect-locIBP}) and then we study the distance between the local densities of two random variables (see Section \ref{sect-dist}). Section \ref{sect-perturbation} is devoted to the results on the perturbation of a Gaussian random variable (see Section \ref{sect-lemma}) and to the study of
the strict positivity and the lower bounds for the density of a general functional on the
Wiener space (see Section \ref{sect-positivity}).
We finally discuss our examples in Section \ref{sect-examples}.

\section{Localized integration by parts formulas}\label{sect-resume}

We consider a probability space $(\Omega ,{\mathscr F},{\mathbb{P}})$ with
an infinite dimensional Brownian motion $W=(W^{n})_{n\in \N_*}$, with $\N_*=\N\setminus\{0\}$, and we use the
Malliavin calculus in order to obtain integration by parts formulas. We
refer to D. Nualart \cite{bib:[N]} for notation and basic results. We denote
by ${\mathbb{D}}^{k,p}$ the space of the random variables which are $k$
times differentiable in Malliavin sense in $L^{p}$ and for a multi-index $%
\alpha =(\alpha _{1},\ldots ,\alpha _{m})\in \N^{m}$ we denote
by $D^{\alpha }F$ the Malliavin derivative of $F$ corresponding to the
multi-index $\alpha .$ So, ${\mathbb{D}}^{m,p}$ is the closure of the space of the simple
functionals with respect to the Malliavin Sobolev norm
$$
\left\Vert F\right\Vert _{m,p}^{p}=\left\Vert F\right\Vert
_{p}^{p}+\sum_{k=1}^{m}
{\mathbb{E}}\big(| D^{(k)}F|^{p}\big)
$$
where
$$
| D^{(k)}F|^2
=\sum_{|\alpha|=k}\int_{[0,\infty)^k}
\left\vert D_{s_{1},\ldots ,s_{k}}^{\alpha}F\right\vert ^{2}\,ds_{1},\ldots ds_{k}.
$$
In the special case $k=1$, we consider the notation
$$
|DF|^2:=|D^{(1)}F|^2
=\sum_{\ell=0}^\infty
\int_{[0,\infty)}\left\vert D^\ell_s F\right\vert ^{2}\,ds,
$$
(for the sake of clearness, we recall that $D^\ell$ stands for the Malliavin derivative w.r.t. $W^\ell$ - and not the derivative of order $\ell$).
Moreover,
for $F=(F^{1},\ldots ,F^{d}),F^{i}\in {\mathbb{D}}^{1,2},$ we let $\sigma _{F}$
denote the Malliavin covariance matrix associated to $F:$%
\begin{equation*}
\sigma _{F}^{i,j}=\langle DF^{i},DF^{j}\rangle
=\sum_{k=1}^{\infty}\int_{0}^{\infty }D_{s}^{k}F^{i}D_{s}^{k}F^{j}ds,\quad
i,j=1,\ldots ,d.
\end{equation*}%
If $\sigma_F$ is invertible, we denote through $\widehat \sigma_F$ the inverse matrix. Finally, as usual, the notation $L$ will be used for the Ornstein-Uhlenbeck operator.

\subsection{Localized representation formulas for the density}\label{sect-locIBP}

Consider a random variable $U$ taking values on $[0,1]$ and set
\begin{equation*}
d\P_{U}=Ud\P.
\end{equation*}%
$\P_{U}$ is a non negative measure (but generally not a probability measure) and we set $\E_U$ the expectation (integral) w.r.t. $\P_U$. For $F\in \D^{k,p}$, we define
$$
\left\Vert F\right\Vert _{p,U}^p =\E_{U}(\left\vert F\right\vert
^{p})\quad \mbox{and} \quad \left\Vert F\right\Vert _{k,p,U}^p =\left\Vert F\right\Vert
_{p,U}^p+\sum_{i=1}^{k}\E_{U}(\vert D^{(i)}F\vert ^{p}).
$$
We assume that $U\in \D^{1 ,\infty }$ and we consider the following condition:
\begin{equation}
m_U(p):=1+\E_{U}(\vert D\ln U\vert ^{p})
<\infty,\quad\mbox{for every $p\in \N$} .
\label{Mall1}
\end{equation}
(\ref{Mall1}) could seem problematic because $U$ may vanish and then
$D(\ln U)$ is not well defined. Nevertheless we make the convention
that $D(\ln U)=\frac{1}{U}DU\,\I_{\{U\neq 0\}}$ (in fact this is the quantity we are really concerned
in). Since $U>0$ $\P_U$-a.s. and $DU$ is well defined, the relation $%
\Vert \ln U\Vert _{1,p,U}<\infty $ makes sense.

We give now the integration by parts formula with respect to $\P_{U}$ (that is, locally) and we study some consequences concerning the regularity of the law starting from the results in Bally and Caramellino \cite{bib:[B.C]} (see also Shigekawa \cite{bib:S} or Malliavin \cite{bib:[M]}). In particular,  for $F\in(\D^{1,\infty})^d$, we will need that the Malliavin covariance matrix $\sigma_F$ is invertible a.s. under $\P_U$, so we call again $\widehat\sigma_F$ the inverse of $\sigma_F$ on the set $\{U\neq 0\}$.

\smallskip

Let $Q_d$ denote the Poisson kernel on $\R^d$: $Q_d$ is  the  fundamental
solution to the equation $\Delta Q_{d}=\delta _{0}$  in $\R^{d}$ ($\delta_0$ denoting the Dirac mass at the origin) and is given by
\begin{equation}\label{den4}
Q_{1}(x)=\max(x,0),\quad
Q_{2}(x)={\cal A}_{2}^{-1}\ln \left\vert x\right\vert \quad \mbox{and}\quad
Q_{d}(x)=-{\cal A}_{d}^{-1}\left\vert x\right\vert ^{2-d},d>2,
\end{equation}%
where for $d\geq 2$, ${\cal A}_{d}$ is the area of the unit sphere in $\R^{d}$.
Then one has

\begin{lemma}\label{1}
 Assume that (\ref{Mall1}) holds. Let $F=(F_{1},...,F_{d})$ be such that $F_{i}\in \D^{2,\infty}$, $i=1,\ldots,d$. Assume that $\det \sigma _{F}>0$ on the set $\{U\neq 0\}$ and moreover
\begin{equation}\label{Law1}
\E_{U}((\det \sigma _{F})^{-p})<\infty \quad \forall p\in \N.
\end{equation}%
Let $\widehat\sigma_{F}$ be the inverse of $\sigma_F$ on the set $\{U\neq 0\}$. Then the following statements hold.

\smallskip

\textbf{A}. For every $f\in C_{b}^{\infty }(R^{d})$ and $V\in \D^{1,\infty }$
one has
\begin{align}
\nonumber
&\E_{U}(\partial_i f(F)\,V)=\E_{U}(f(F)H_{i,U}(F,V)),\quad i=1,\ldots,d,\mbox{with }\\
&H_{i,U}(F,V)=\sum_{j=1}^d\Big(V\widehat \sigma^{ji} _{F}LF^j-\left\langle D(V\widehat\sigma^{ji} _{F}),DF^j\right\rangle
-V\widehat\sigma^{ji}_{F}\left\langle D\ln U,DF^j\right\rangle\Big).  \label{Mall5}
\end{align}%

\smallskip

\textbf{B}. Let $Q_d $ be the Poisson kernel in $\R^{d}$ given in (\ref{den4}). Then for every $p>d$ one has
\begin{equation}
\E_{U}(\left\vert \nabla Q_d (F-x)\right\vert ^{\frac{p}{p-1}})^{\frac{p-1}{%
p}}\leq C_{p,d}\E_{U}(\left\vert H_{U}(F,1)\right\vert ^{p})^{k_{p,d}}
\label{Mall5'}
\end{equation}%
where $C_{p,d}$ is a universal constant depending on $p$ and $d$ and
$k_{p,d}=(d-1)/(1-d/p)$.

\smallskip

\textbf{C}. Under $\P_U$, the law of $F$ is absolutely continuous and has a
continuous density $p_{F,U}$ which may be represented as%
\begin{equation}
p_{F,U}(x)=\sum_{i=1}^d\E_{U}(\partial_i Q_d (F-x)H_{i,U}(F,1)).
\label{Mall5''}
\end{equation}%
Moreover, for every $p>d$ there exist constants $C>0$ and $q>1$ depending on $p,d$ such that%
\begin{equation}
p_{F,U}(x)\leq C\si _{F,U}(p)^{q}n_{F,U}(p)^{q}m_U(p)^q  \label{Bis1}
\end{equation}%
with $m_U(p)$ given in (\ref{Mall1}),
\begin{equation}\label{F,U}
\si _{F,U}(p)=1+\E_{U}(\left\vert \det \sigma_{F}\right\vert ^{-p})\quad
and\quad n_{F,U}(p)=1+\left\Vert F\right\Vert _{2,p,U}+\|LF\|_{p,U}.
\end{equation}
Finally, if $V\in \D^{1,\infty }$ then
\begin{equation}
p_{F,UV}(x)\leq C\si _{F,U}(p)^{q}n_{F,U}(p)^{q}m_U(p)^q\left\Vert V\right\Vert
_{1,p,U},  \label{Bis3}
\end{equation}
in which $p>d$ and $C>0$, $q>1$ are suitable constants depending on $p,d$.
\end{lemma}

\textbf{Proof}. \textbf{A.} The standard integration by parts formula in Malliavin
calculus gives (vector notations)
\begin{equation*}
\E_{U}(\nabla f(F)V)=\E(\nabla f(F)UV)=\E(f(F)H(F,UV))
\end{equation*}%
where, setting $DU=U\times D(\ln U)$, one has
\begin{eqnarray*}
H(F,UV) &=&VU\widehat\sigma _{F}LF-\left\langle D(VU\widehat\sigma _{F}),DF\right\rangle \\
&=&U(V\widehat\sigma _{F}LF-\left\langle D(V\widehat\sigma _{F}),DF\right\rangle )-V\widehat\sigma
_{F}\left\langle D\ln U,DF\right\rangle ),
\end{eqnarray*}%
So, $H(F,UV)=UH_{U}(F,V)$, and (\ref{Mall5}) is proved.

\smallskip

\textbf{B.} This point straightforwardly follows from the results and the techniques  in Bally and Caramellino \cite{bib:[B.C]}.

\smallskip

\textbf{C.} (\ref{Mall5''}) again follows from \cite{bib:[B.C]}, while (\ref{Bis1}) is a consequence of the inequality
\begin{equation}
\|H_{U}(F,V)\|_{p,U}\leq C\si
_{F,U}(p)^{q}n_{F,U}(p)^{q}m_U(p)\left\Vert V\right\Vert _{1,p,U},  \label{Bis4'}
\end{equation}
holding for suitable $C>0$ and $p,q>1$ depending on $d$ only. This can be proved by applying the H\"older inequality to (\ref{Mall5}) (further details can be found in the proof of next Proposition \ref{prop-dist}). So, by using the H\"older inequality to (\ref{Mall5''}) and by considering both (\ref{Mall5'}) and (\ref{Bis4'}), one gets (\ref{Bis1}).
Finally, in order to prove (\ref{Bis3}) we formally write (the rigorous arguments can be found in \cite{bib:[B.C]})
\begin{align*}
p_{F,UV}(x) =&\E_{UV}(\delta _{0}(F-x))=\E_{UV}(\triangle Q_d(F-x))=\E_{U}(\triangle Q_d (F-x)V) \\
=&\E_{U}(\left\langle \nabla Q_d (F-x),H_{U}(F,V\right\rangle ).
\end{align*}%
Then using (\ref{Mall5'}) and(\ref{Bis4'}) one obtains (\ref{Bis3}). $\square
$

\subsection{The distance between two density functions}\label{sect-dist}

We compare now the densities of the laws of two random variables under $\P_{U}.$

\begin{proposition}\label{prop-dist}
 Assume that (\ref{Mall1}) holds. Let $F=(F_{1},...,F_{d})$ and $%
G=(G_{1},...,$ $G_{d})$ be such that $F_{i},G_{i}\in \D^{2,\infty}$, $i=1,\ldots,d$, and
\begin{equation*}
\si _{F,G,U}(p):=1+\sup_{0\leq \varepsilon \leq 1}\E_{U}((\det \sigma
_{G+\varepsilon (F-G)})^{-p}))<\infty,\quad\forall p\in\N.
\end{equation*}%
Then under $\P_U$ the laws of $F$ and $G$
are absolutely continuous with respect to the Lebesgue measure with
continuous densities $p_{F,U}$ and $p_{G,U}$ respectively. Moreover, for every $p>d$ there exist two constants $C>0$, $q>1$ depending on $p,d$ such that
\begin{equation}\label{Mall2}
\left\vert p_{F,U}(y)-p_{G,U}(y)\right\vert \leq C\,\si
_{F,G,U}(p)^{q}n_{F,G,U}(p)^{q}m_U(p)^q\,\|\Delta_2(F,G)\|_{p,U}
\end{equation}%
with $m_U(p)$ given in (\ref{Mall1}) and
\begin{equation} \label{Mall2'}
\begin{array}{l}
\displaystyle
\Delta _{2}(F,G) =\left\vert D(F-G)\right\vert +\left\vert
D^{(2)}(F-G)\right\vert +\left\vert L(F-G)\right\vert ,  \smallskip\\
\displaystyle
n_{F,G,U}(p) =1+\left\Vert F\right\Vert _{2,p,U}+\left\Vert G\right\Vert
_{2,p,U}+\|LF\|_{p,U}+\|LG\|_{p,U}.
\end{array}
\end{equation}
Finally, recalling that $\left\vert U\right\vert \leq 1$ almost surely, using Meyer's
inequality one has
\begin{equation}
\left\vert p_{F,U}(y)-p_{G,U}(y)\right\vert \leq C\si
_{F,G,U}(p)^{q}m_U({p})^q(1+\left\Vert F\right\Vert _{2,p}+\left\Vert
G\right\Vert _{2,p})^{q} \left\Vert F-G\right\Vert _{2,p},
\label{Mall2''}
\end{equation}
for $p>d$ and $C>0$, $q>1$ depending on $p,d$.
\end{proposition}

\textbf{Proof}. Throughout this proof, $C ,p,q$ will denote constants that can vary from line to line.

\smallskip

By applying Lemma \ref{1}, we first notice that under $\P_U$ the laws of $F$ and $G$
are both absolutely continuous with respect to the
Lebesgue measure and the densities can be written as
\begin{equation}\label{pFGU}
\begin{array}{l}
p_{F,U}(y)=\E_{U}(\left\langle \nabla Q_d(F-y),H_{U}(F,1)\right\rangle )\mbox{ and }\smallskip\\
p_{G,U}(y)=\E_{U}(\left\langle \nabla Q_d(G-y),H_{U}(G,1)\right\rangle ).
\end{array}
\end{equation}%

\textbf{Step 1}. We prove that for $V\in \D^{1,\infty }$, on the set $\{U\neq 0\}$ one has
\begin{equation} \label{Mall3}
\left\vert H_{U}(F,V)-H_{U}(G,V)\right\vert \leq C\, A_{F,G}\,
B_{F,G}\, (1+\left\vert D\ln U\right\vert )(\left\vert V\right\vert
+\left\vert DV\right\vert )\times \Delta _{2}(F,G)
\end{equation}%
where on the set $\{U\neq 0\}$ (that is, where the inverse Malliavin covariance matrices
$\widehat\sigma_F$ and $\widehat\sigma_G$ are actually well defined) the above quantities are equal to
\begin{align*}
&A_{F,G} =(1\vee \det \widehat\sigma
_{F})^2(1\vee\det \widehat\sigma _{G})^{2},\\
&B_{F,G} =(1+\left\vert DF\right\vert +\left\vert DG\right\vert +\left\vert
D^{(2)}F\right\vert +\left\vert D^{(2)}G\right\vert )^{\alpha_d}(1+\left\vert
LF\right\vert +\left\vert LG\right\vert ),
\end{align*}
$\alpha_d$ denoting a positive constant depending on $d$ only.
So, we work on the set $\{U\neq 0\}$. We first notice that
\begin{equation}\label{1bisbis}
\left\vert \widehat\sigma _{F}^{i,j}-\widehat\sigma _{G}^{i,j}\right\vert
\leq
C(1\vee \det \widehat\sigma
_{F})(1\vee\det \widehat\sigma _{G})\left\vert D(F-G)\right\vert(1+\left\vert DF\right\vert +\left\vert
DG\right\vert )^{2d-1}
\end{equation}
and moreover,
\begin{equation}\label{2bisbis}
\begin{array}{c}
\left\vert D\widehat\sigma _{F}^{i,j}-D\widehat\sigma _{G}^{i,j}\right\vert
\leq  C(1\vee \det \widehat\sigma
_{F})^2(1\vee\det \widehat\sigma _{G})^{2} (\left\vert D(F-G)\right\vert +\left\vert D^{(2)}(F-G)\right\vert )\times\smallskip\\
\times (1+\left\vert DF\right\vert +\left\vert
DG\right\vert +\left\vert D^{(2)}F\right\vert +\left\vert D^{(2)}G\right\vert
)^{6d-3}
\end{array}
\end{equation}
The proof of (\ref{1bisbis}) and (\ref{2bisbis}) is technical but standard, so we postpone it to Appendix \ref{bisbis}. By  a straightforward computation, including the use of (\ref{1bisbis}) and (\ref{2bisbis}),  one gets (\ref{Mall3}).
Hence, (\ref{Mall3}) and the H\"{o}lder inequality give
\begin{equation}
\|H_{U}(F,V)-H_{U}(G,V)\|_{p,U}\leq
Cn_{F,G,U}(p^{\prime })^{q}m_U(p^{\prime })\left\Vert V\right\Vert
_{1,p^{\prime },U}\|\Delta_2(F,G)\|_{p',U}
\label{Mall6}
\end{equation}%
for some $p'>p$.

\smallskip

\textbf{Step 2.} By using arguments similar to the ones developed in Step 1, we get
\begin{equation}
\|H_{U}(F,V)\|_{p,U}\leq C\,\si
_{F,U}(p')^{q}n_{F,U}(p')^{q}m_U(p')\left\Vert V\right\Vert _{1,p',U},  \label{Bis4}
\end{equation}
$n_{F,U}(p')$ and $\si
_{F,U}(p')$ being defined in (\ref{F,U}).
So, by taking $p>d$ 
in (\ref{Mall5'}) and by using (\ref{Bis4}) with $V=1$ one gets
\begin{equation}
\|\nabla Q_{d}(F-y)\|_{p/(p-1),U}\leq C\,
\si _{F,U}(p')^{q}\,n_{F,U}(p')^{q}m_U(p')^q,  \label{Mall7}
\end{equation}
with $p'>p>d$.

\smallskip

\textbf{Step 3}. By using (\ref{pFGU}), we can write
\begin{align*}
p_{F,U}(y) -p_{G,U}(y)=&\E_{U}(\langle \nabla Q_{d}(F-y)-\nabla
Q_{d}(G-y),H_{U}(G,1)\rangle ) +\\
&+\E_{U}(\langle\nabla Q_{d}(F-y),H_{U}(F,1)-H_{U}(G,1)\rangle) \\
=&:I+J.
\end{align*}%
Using (\ref{Mall6}) we obtain%
\begin{equation*}
|J|\leq C\,\si _{F,G,U}(p)^{q} n_{F,G,U}(p)^{q}
m_U(p)^{q}\|\Delta _{2}(F,G)\|_{p,U},
\end{equation*}%
with $p>d$ and $C>0,q>1$ depending on $p$ and $d$. We study now the quantity $I$. For $\lambda \in \lbrack 0,1]$ we denote $F_{\lambda }=G+\lambda (F-G)$ and
we use Taylor's expansion to obtain%
\begin{equation*}
I=\sum_{k,j=1}^{d}R_{k,j}\quad \mbox{with}\quad R_{k,j}=\int_{0}^{1}\E_{U}(\partial
_{k}\partial _{j}Q_{d}(F_{\lambda }-y)H_{j,U}(G,1)(F-G)_{k})d\lambda .
\end{equation*}%
Let $V_{k,j}=H_{j,U}(G,1)(F-G)_{k}.$ Using again the integration by parts
formula (with respect to $F_{\lambda })$ we obtain%
\begin{equation*}
R_{k,j}=\int_{0}^{1}\E_{U}\big(\partial _{j}Q_{d}(F_{\lambda
}-y)H_{k,U}(F_{\lambda },V_{k,j})\big)d\lambda .
\end{equation*}
Now, one has $\E_U((\det \sigma_{F_\lambda})^{-p})\leq \si_{F,G,U}(p)<\infty$ for every $\lambda\in [0,1]$ and $p\geq 1$. So, we can use (\ref{Mall7}) and (\ref{Bis4}) with $F=F_\lambda$, and we get
\begin{eqnarray*}
\left\vert R_{k,j}\right\vert &\leq &C\si
_{U,F,G}(p)^{q}n_{U,F,G}(p)^{q}m_U(p)^q\left\Vert V_{k,j}\right\Vert _{1,p,U} \\
&\leq &C^{\prime }\si _{U,F,G}(p')^{q'}n_{U,F,G}(p^{\prime
})^{q'}m_U(p')^{q'}\|\Delta_2(F,G)\|_{p',U},
\end{eqnarray*}%
with $p'>p>d$ and $C'>0,q'>1$ depending on $p', d$. The statement now easily follows.
$\square $

\medskip

\begin{example}\label{U-psi}
We give here an example of localizing function giving rise to a localizing random variable $\bar U$ that satisfies (\ref{Mall1}). For $a>0$, set $\psi
_{a}:{\mathbb{R}}\rightarrow {\mathbb{R}}_{+}$ as
\begin{equation}\label{Mall10}
\psi _{a}(x)=1_{|x|\leq a}+\exp \Big(1-\frac{a^2}{a^2-(x-a)^2}\Big)1_{a<|x|< 2a}.
\end{equation}%
Then $\psi _{a}\in C^1_b(\R)$, $0\leq \psi_a\leq 1$  and  for every $p\geq 1$ one has
$$
\sup_{x}|(\ln\psi_a(x))'|^p\psi_a(x)\leq \frac {4^p}{a^{p}}\,\sup_{t\geq 0}(t^{2p}e^{1-t})<\infty.
$$
For  $\Theta _{i}\in \D^{1,\infty}$ and $a_{i}>0$, $i=1,...,\ell$, we define
\begin{equation}
\bar U=\prod_{i=1}^{\ell}\psi _{a_{i}}(\Theta _{i}).  \label{Mall10'}
\end{equation}%
Then $\bar U\in\D^{1,\infty}$, $\bar U\in [0,1]$ and  (\ref{Mall1}) holds. In fact, one has
\begin{align*}
|D\ln \bar U|^p\bar U
=&\Big|\sum_{i=1}^\ell(\ln\psi_{a_i})'(\Theta_i)D\Theta_i\Big|^p\prod_{j=1}^\ell\psi_{a_j}(\Theta_j)\\
\leq&\Big(\sum_{i=1}^\ell|(\ln\psi_{a_i})'(\Theta_i)|^2\Big)^{p/2}\Big(\sum_{i=1}^\ell|D\Theta_i|^2\Big)^{p/2}
\prod_{j=1}^\ell\psi_{a_j}(\Theta_j)\\
\leq& c_p \sum_{i=1}^\ell|(\ln\psi_{a_i})'(\Theta_i)|^p\psi_{a_i}(\Theta_i)\times |D\Theta|^p\\
\leq& C_{p}\sum_{i=1}^\ell \frac 1{a_i^{p}}\,|D\Theta|^p
\end{align*}
for a suitable $C_{p}>0$, so that
\begin{equation}\label{Mall11}
\E(|D\ln \bar U|^p\bar U)\leq C_{p}\sum_{i=1}^\ell \frac 1{a_i^{p}}\times \E(|D\Theta|^p)
\leq C_{p}\sum_{i=1}^\ell \frac 1{a_i^{p}}\times\|\Theta\|_{1,p}^p<\infty.
\end{equation}
\end{example}

Using the localizing function in (\ref{Mall10}) and by applying Proposition \ref{prop-dist} we get the following result.

\begin{theorem}\label{th-dist}
Assume that (\ref{Mall1}) holds. Let $%
F=(F_{1},...,F_{d})$ and $G=(G_{1},...,G_{d})$\ with $F_{i},G_{i}\in \D^{2,\infty}$
and such that for every $p\in\N$ one has
\begin{equation*}
\si _{F,U}(p):=1+\E_{U}((\det \sigma _{F})^{-p}))<\infty \quad and\quad
\si _{G,U}(p):=1+\E_{U}((\det \sigma _{G})^{-p}))<\infty.
\end{equation*}%
Then under $\P_U$, the laws of $F$ and $G$ are absolutely continuous with respect to the Lebesgue measure,
with continuous densities $p_{F,U}$ and $p_{G,U}$ respectively, and for $p>d$ there exist two constant $C>0$
and $q>1$ depending on $p,d$ such that
\begin{equation} \label{Mall12}
\left\vert p_{F,U}(y)-p_{G,U}(y)\right\vert \leq C\,(\si _{G,U}(p)\vee
\si _{F,U}(p))^{q} n_{F,G,U}(p)^{q}m_U(p)^{q}\times
\|\Delta_2(F,G)\|_{p,U}
%
\end{equation}
with $n_{F,G,U}(p)$ and $\Delta _{2}(F,G)$ given in (\ref{Mall2'}) and $%
m_U(p)$ given in (\ref{Mall1}).

\end{theorem}

\textbf{Proof}. Set $R=F-G.$ By using (\ref{b}) (see Appendix \ref{bisbis}), for every $\lambda\in [0,1]$ one gets
\begin{equation}\label{alpha-d}
\det
\sigma _{G+\lambda R}\geq \det \sigma _{G}-\alpha_d\left\vert DR\right\vert
\left\vert DG\right\vert (1+\left\vert DF\right\vert +\left\vert
DG\right\vert )^{2d-1},
\end{equation}
for a suitable $\alpha_d>0$ depending on $d$ only.
For $\psi_a$ as in (\ref{Mall10}), we define
\begin{equation*}
V=\psi _{1/4}(H)\quad \mbox{with}\quad H=\frac{\alpha_d\left\vert DR\right\vert
\left\vert DG\right\vert (1+\left\vert DF\right\vert +\left\vert
DG\right\vert )^{2d-1}}{\det \sigma _{G}}
\end{equation*}%
so that if $V\neq 0$ then $\det \sigma _{G+\lambda R}\geq \frac{1}{2}\det
\sigma _{G}$. It follows that $\si _{F,G,UV}(p)\leq C\si _{G,U}(p)$, $C$ denoting a suitable positive constant (which will vary in the following lines).
We also have $m_{UV}(p)\leq C(m_U(p)+\E(UV|D\ln V|^{p}))$ and by (\ref{Mall11})\ we have%
\begin{equation*}
\E(UV|D\ln V|^{p})\leq C \left\Vert DH\right\Vert
_{p,U}^p\leq C\, n_{F,G,U}(\bar p)^{\bar q} \si _{G,U}(\bar p)^{\bar q}
\end{equation*}%
for some $\bar p, \bar q$, so that
$m_{UV}(p)\leq C\,m_U(p)n_{F,G,U}(\bar p)^{\bar q} \si _{G,U}(\bar p)^{\bar q}$.
So, we can apply (\ref{Mall2}) with localization $UV$ and we get
\begin{equation}
\left\vert p_{F,UV}(y)-p_{G,UV}(y)\right\vert \leq C\si _{G,U}(p)^{q} n_{F,G,U}(p)^{q} m_U(p)^{q}\|\Delta
_{2}(F,G)\|_{p,U},
\label{Mall17}
\end{equation}%
with $p>d$ and $C>0,q>1$ depending on $p,d$.
 We write now%
$$
\left\vert p_{F,U}(y)-p_{G,U}(y)\right\vert \leq \left\vert
p_{F,UV}(y)-p_{G,UV}(y)\right\vert +\left\vert p_{F,U(1-V)}(y)\right\vert
+\left\vert p_{G,U(1-V)}(y)\right\vert ,
$$
and we have already seen that the first addendum on the r.h.s. behaves as desired. So, it suffices to show that also the remaining two terms have the right behavior. To this purpose, we use (\ref{Bis3}). We have
\begin{equation*}
\left\vert p_{F,U(1-V)}(y)\right\vert \leq \si _{F,U}(p)^q
n_{F,1,U}(p)^qm_U(p)^q\times \left\Vert 1-V\right\Vert _{1,p,U}.
\end{equation*}%
We recall that $1-V\neq 0$ implies that $H\geq 1/8$, so that
\begin{align}
\|1-V\|_{1,p,U}^p
&=\E_U(|1-V|^p)+\E_U(|DV|^p)
\leq C\big(\P_U(H>1/8)+\E_U(V|D\ln V|^p)\big)\nonumber\\
&\leq C\big(\E_U(H^p)+\E_U(|DH|^p)\big)\label{1-V}
\end{align}
in which we have used (\ref{Mall11}). Now, one has
\begin{align*}
\E_U(|H|^p)
&\leq C\,\gamma_{G,U}(\bar p)^{\bar q}n_{F,G,U}(\bar p)^{\bar q}\,\E_U(|D(F-G)|^{2p})^{1/2}\quad\mbox{and}\\
\E_U(|DH|^p)
&\leq C\,\gamma_{G,U}(\bar p)^{\bar q}n_{F,G,U}(\bar p)^{\bar q}\,\big(\E_U(|D(F-G)|^{2p})^{1/2}+\E_U(|D^{(2)}(F-G)|^{2p})^{1/2}\big)
\end{align*}
and by inserting above we get
$$
\|1-V\|_{1,p,U}
\leq C\,\gamma_{G,U}(\bar p)^{\bar q}n_{F,G,U}(\bar p)^{\bar q}\,\|\Delta_2(F,G)\|_{2p,U}.
$$
This gives
\begin{equation*}
\left\vert p_{F,U(1-V)}(y)\right\vert \leq C\big(\si _{F,U}(p)\vee\si _{G,U}(p)\big)^q n_{F,G,U}(p)^qm_U(p)^q \|\Delta_2(F,G)\|_{p,U}
\end{equation*}
for $p>d$ and suitable constants  $C>0$ and $q>1$ depending on $p,d$. And similarly we get
\begin{equation}\label{Mall17bis}
\left\vert p_{G,U(1-V)}(y)\right\vert \leq C\,\si _{G,U}(p)^q n_{F,G,U}(p)^qm_U(p)^q \|\Delta_2(F,G)\|_{p,U},
\end{equation}
with the same constraints for $p,C,q$. The statement now follows. $\square$

\medskip

An immediate consequence of the above proof consists in a lower bound for $p_{F,U}$ that does not involve $\si_{F,U}(p)$, that is the Malliavin covariance matrix of $F$. This can be done thanks to the localizing r.v. $V$ underlying the proof of Theorem \ref{th-dist}. In fact, one has

\begin{proposition}\label{Mall-rem}
Under the hypotheses of Theorem \ref{th-dist}, for every $p>d$ there exist two constants $C>0$ and $q>1$
depending on $p,d$ such that
$$
p_{F,U}(y)\geq p_{G,U}(y)-
C\,\si _{G,U}(p)^{q} n_{F,G,U}(p)^{q}m_U(p)^{q}\times
\|\Delta_2(F,G)\|_{p,U}.
$$
\end{proposition}

\textbf{Proof.}
Let $V$ be the localizing r.v. as in the proof of Theorem \ref{th-dist}. We can write
\begin{align*}
p_{F,U}(y)
&\geq p_{F,UV}(y)\geq p_{G,UV}(y)-|p_{F,UV}(y)-p_{G,UV}(y)|\\
&= p_{G,U}(y)-p_{G,U(1-V)}(y)-|p_{F,UV}(y)-p_{G,UV}(y)|.
\end{align*}
The statement now follows from (\ref{Mall17bis}) and (\ref{Mall17}). $\square$

\section{Small perturbations of a Gaussian random variable}\label{sect-perturbation}

\subsection{Preliminary estimates}\label{sect-lemma}
We consider here a r.v. of the type $F=x+G+R\in {\mathbb{R}}^{d}$ where $R\in
{\mathbb{D}}^{2,\infty }$ and
\begin{equation*}
G=\sum_{j=1}^{\infty }\int_{0}^{\infty }h_{j}(s)dW_{s}^{j}
\end{equation*}%
with $h_{j}\,:\,[0,+\infty)\to {\mathbb{R}}^d$ deterministic and
square integrable. Then $G$ is a
centered Gaussian random variable of covariance matrix $M_{G}=(M_{G}^{k,p})_{k,p=1,\ldots,d}$, with
\begin{equation*}
M_{G}^{k,p}=\int_{0}^{\infty }\langle h^{k}(s),h^{p}(s)\rangle
ds=\sum_{j=1}^{\infty }\int_{0}^{\infty }h_{j}^{k}(s)h_{j}^{p}(s)ds, \quad k,p=1,\ldots,d.
\end{equation*}
We assume that $M_{G}$ is invertible and we denote by $g_{M_{G}}$ the
density of $G$ that is%
\begin{equation*}
g_{M_{G}}(y)=\frac{1}{(2\pi )^{d/2}\sqrt{\det M_{G}}}\exp (-\left\langle
M_{G}^{-1}y,y\right\rangle ).
\end{equation*}%
Our aim is to give estimates of the  density of $F$ in terms of $g_{M_{G}}$. To this purpose, we use a localizing r.v. $U$ of the form (\ref{Mall10'}).

\begin{proposition}\label{exp-est}
Let $\psi_a$ be the function in (\ref{Mall10}) and set
\begin{equation}\label{Law5'}
U=\psi_{c_*^2/2}(\vert D\overline{R}\vert ^{2})\quad\mbox{with}\quad
 \overline{R}=M_G^{-1/2}R.
\end{equation}
where $c_*$ is such that $\alpha_d(1+2d)^{d}c_*(1+c_*)^{d-1}\leq 1/2$, $\alpha_d$ denoting the constant in (\ref{alpha-d}).
Then the following statements hold.

\smallskip

i) Under $\P_U$, the law of $F=x+G+R$ has a smooth density $p_{F,U}$ and one has
\begin{equation*}
\sup_{y\in \R^{d}}\left\vert p_{F,U}(y)-g_{M_{G}}(y-x)\right\vert
\leq \varepsilon (M_{G},R),\quad y\in\R^d,
\end{equation*}%
where%
\begin{equation*}
\varepsilon (M_{G},R):=\frac{c_{d}}{\sqrt{\det M_{G}}}(1+\left\Vert
\overline{R}\right\Vert _{2,q_{d}})^{\ell _{d}}\left\Vert \overline{R}%
\right\Vert _{2,q_{d}}.
\end{equation*}%
Here $c_{d}>0$ and $q_{d},\ell _{d}>1$ are universal constants depending on $d$ only.

\smallskip

ii) If the law of $F$ under ${\mathbb{P}}$ has a density $p_{F}$,
then one has
\begin{equation*}
p_{F}(y)\geq g_{M_{G}}(y-x)-\varepsilon (M_{G},R),\quad y\in
\R^{d}.
\end{equation*}%
\end{proposition}

\textbf{Proof}. $i)$
Suppose first that $x=0$ and  $M_G=I$, $I$ denoting the identity matrix, so that $\overline{R}=R$. We notice that $\det \sigma_G=1$, which gives  $\gamma_{G,U}(p)\leq 2$ for every $p$, and $|DG|^2=d$. Moreover,
on the set $\{U\neq 0\}$ one has $|DR|\leq c_*$ and by using (\ref{alpha-d}) straightforward computations give
$$
\det\sigma_F\geq 1-\alpha_d(1+2d)^{d}|DR|(1+|DR|)^{d-1}
\geq 1- \alpha_d(1+2d)^{d}c_*(1+c_*)^{d-1}\geq \frac 12.
$$
It then follows that $\gamma_{F,U}(p)\leq 1+2^{p}<\infty$ for every $p$. Moreover, by (\ref{Mall11}) one has
$m_U(p)\leq 1+\|R\|_{2,2p}^p$. We can then apply Theorem \ref{th-dist} to the pair $F$ and $G$, with
localizing r.v. $U$. By straightforward computations and the use of the Meyer inequality, one has
$n_{F,G,U}(p)\leq C_p(1+\|R\|_{2,2p})^{2p}$ and $\|\Delta _{2}(F,G))\|_p \leq C_p\|R\|_{2,2p}$, with
$n_{F,G,U}(p)$ and $\Delta _{2}(F,G)$ given in (\ref{Mall2'}). Therefore, by applying (\ref{Mall12}) with
$p=d+1$, one has
$$
\left\vert p_{F,U}(y)-p_{G,U}(y)\right\vert \leq
\bar c_1\,
(1+\|R\|_{2,\bar q_1})^{\bar \ell_1}\,
\|R\|_{2,\bar q_1}\quad\mbox{for every $y\in\R^d$},
$$
where $\bar c_1>0$ and $\bar q_1,\bar \ell_1>1$ are constants depending on $d$ only.
It remains to  compare $p_{G,U}$ with $p_{G}=g_I$: from (\ref{Bis3}) (applied with $U=1$, $V=1-U$ and $F=G$) we immediately have
$$
|p_{G,U}(x)-p_{G}(x)|
=p_{G,1-U}(x)\leq C\si _{G,1}(p)^{q}n_{G,1}(p)^{q}m_1(p)^q\left\Vert 1-U\right\Vert
_{1,p}
<C\|1-U\|_{1,p},
$$
for $p>d$ and $C>, q>1$ depending on $p,d$ only. Now, recalling that $U\neq 1$ for $|DR|>c_*/\sqrt 2$, as
already seen in the proof of Theorem \ref{th-dist} (see (\ref{1-V})) we have
$$
\|1-U\|_{1,p}^p
\leq C\big(\E(|DR|^{2p})+\E(|D|DR|^{2}|^p)\big),
$$
so that taking $p=d+1$
$$
|p_{G,U}(y)-p_{G}(y)|\leq \bar c_2
(1+\|R\|_{2,\bar q_2})^{\bar \ell_2}\,
\|R\|_{2,\bar q_1}\quad\mbox{for every $y\in\R^d$},
$$
with $\bar c_2>0$ and $\bar q_1,\bar q_2,\bar \ell_2>1$ depending on $d$ only, and the statement follows. As
for the general case, it suffices to apply the already proved estimate to $\overline{F}=M_G^{-1/2}(F-x)$,
$\overline{G}=M_G^{-1/2}G$ and $\overline{R}=M_G^{-1/2}R$ and then to use the change of variable theorem.

\smallskip

$ii)$ It immediately follows from $p_F(y)\geq p_{F,U}(y)$.
$\square$

 \medskip

\subsection{Main results}\label{sect-positivity}

In this section, we consider a time interval of the type $[T-\delta,T]$, where
$T>0$ is a fixed horizon and $0<\delta\leq  T$, and
we use the Malliavin calculus with respect to
$W_{s},s\in
[ T-\delta ,T]$. In particular, we take conditional expectations with respect to $%
\mathcal{F}_{T-\delta }.$ Therefore, for $V=(V^{1},...,V^{d})$, $V_i\in \D^{L,p}$,
we define the following conditional
Malliavin  Sobolev norms:
\begin{equation}\label{cond-mall-norm}
\left\Vert V\right\Vert _{\delta ,L,p}^p =\E(\left\vert
V\right\vert ^{p}\mid \mathcal{F}_{T-\delta })
+\sum_{l=1}^{L}\E\big(\vert
D^{(l)}V\vert^{p}
\,\big |\,\mathcal{F}_{T-\delta }\big).
\end{equation}%

Let $F$ denote a $d$-dimensional functional  on the Wiener space which is measurable w.r.t. $\mathcal{F}_{T}$ and assume that for $\delta
\in (0, T]$ the following decomposition holds:
\begin{equation}\label{deco}
F=F_{T-\delta }+G_{\delta }+R_{\delta }
\end{equation}%
where $F_{T-\delta }$ is measurable w.r.t. $\mathcal{F}_{T-\delta }$, $R_{\delta }\in (\D^{2,\infty })^{d}$
and
$$
G_{\delta }=\sum_{k=1}^{\infty }\int_{T-\delta }^{T}h_{\delta
}^{k}(s)dW_{s}^{k}.
$$
Here $h_{\delta }^{k}(s),s\in \lbrack
T-\delta ,T]$ are progressively
measurable processes such that $h_{\delta }^{k}(s)$ is $\mathcal{F}%
_{T-\delta }$-measurable for every $s\in[T-\delta,T]$ and $\sum_{k=1}^{\infty }\int_{T-\delta
}^{T}\vert h_{\delta }^{k}(s)\vert ^{2}ds<\infty $
a.s.
In particular, conditionally on
$\mathcal{F}_{T-\delta },$ the random variable $G_{\delta }$ is
centered and Gaussian with covariance matrix
\begin{equation*}
C_{\delta }^{ij}=\sum_{k=1}^{\infty }\int_{T-\delta }^{T}h_{\delta
}^{k,i}(s)h_{\delta }^{k,j}(s)ds\quad 1\leq i,j\leq d.
\end{equation*}

On the set $\{\det C_\delta\neq 0\}\in\cl F_{T-\delta}$, we define the (random) norm
\begin{equation*}
\vert x\vert _{\delta }:=\vert C_{\delta
}^{-1/2}x\vert,\quad x\in \R^d
\end{equation*}%
and for $q\in\N$,  we consider the following (random) quantity
\begin{equation}\label{theta}
\theta_{\delta,q}=\|C_\delta^{-1/2}R_\delta\|_{\delta,2,q}.
\end{equation}
Set now $\overline{\P}_\delta(\omega,\cdot)$ the measure induced
by $\overline{\E}_\delta(\omega,X)=\E(X\psi(|DR_\delta|^2)\,|\,\cl F_{T-\delta})(\omega)$,
where $\psi=\psi_{1/8}$ is as in (\ref{Law5'}).
By developing in a conditional form the arguments as in the proof of Proposition \ref{exp-est},
on the set $\{\det C_\delta\neq 0\}$ one gets that
under $\overline{\P}_\delta(\omega,\cdot)$
the law of $F$ has a regular density w.r.t. the Lebesgue measure.
Therefore, there exists a function
$\bar p_{F,\delta}(\omega,z)$ which is regular as a function of
$z$ and such that
\begin{equation}\label{pbardelta0}
\E\big(f(F)\psi(|DR_\delta|^2)\,|\,\cl F_{T-\delta})(\omega)
=\int f(z)\bar p_{F,\delta}(\omega,z)dz,\quad \omega\in\{\det C_\delta\neq 0\}
\end{equation}
for any measurable and bounded function $f$.

\smallskip

Now, let us introduce  the following sets:
for $y\in \R^d$ and $r>0$,
we define
\begin{align}
&\Gamma_{\delta,r}(y)=\big\{\vert F_{T-\delta }-y\vert _{\delta }\leq
r\big\}\cap \big\{\det C_\delta\neq 0\big\}\cap \big\{\theta _{\delta ,q_d}\leq
a_d\,e^{-r^2}\big\}\label{gamma0}\\
&
\widetilde\Gamma_{\delta,r}(y)=\big\{\vert F_{T-\delta }-y\vert _{\delta }\leq
r/2\big\}\cap \big\{\det C_\delta\neq 0\big\}\cap \big\{\theta _{\delta ,q_d}\leq
a_d\,e^{-r^2}\big\},\label{gamma}
\end{align}
where
$$
a_d
=\frac 1{c_d2^{\ell_d+1}(2\pi)^{d/2}}
$$
and $q_d$, $\ell_d$ and $c_d$ are the universal constants defined in $(i)$ of Proposition \ref{exp-est}.
Then we have

\begin{lemma}\label{lemma-pos}
For $\delta\in (0,T]$, let decomposition (\ref{deco}) hold and for $y\in\R^d$, $r>0$, let
$\Gamma_{\delta,r}(y)$ be the set in (\ref{gamma0}). Then  for every non negative and measurable function
$f\,:\,\R^d\to\R$ and $\omega\in\{\det C_\delta\neq 0\}$ one has
$$
\E(f(F)\,|\,\cl F_{T-\delta})(\omega)
\geq \frac{e^{-r^2}}{2(2\pi)^{d/2}}\,(\det C_\delta)^{-1/2}\int f(z)\I_{\Gamma_{\delta,r}(z)}dz.
$$
\end{lemma}

\textbf{Proof}. Let $\omega \in \{\det C_{\delta }\neq 0\}.$
By using (\ref{pbardelta0}),
for any measurable and non negative function $f$ we have
\begin{align*}
\E(f(F)\,|\, \mathcal{F}_{T-\delta })(\omega )
&\geq \E(f(F)\psi (\vert
DR_{\delta }\vert ^{2})\mid \mathcal{F}_{T-\delta })(\omega )
=\int f(z)\overline{p}_{F,\delta }(\omega ,z)dz \\
&\geq \int f(z)\overline{p}_{F,\delta }(\omega ,z)\I_{\Gamma _{\delta
,r}(z)}dz
\end{align*}
Using Proposition \ref{exp-est}  in a conditional form (with respect to $\mathcal{F}%
_{T-\delta })$ we obtain
\begin{equation*}
\overline{p}_{F,\delta }(\omega ,z)\geq g_{C_{\delta }(\omega)}(z-F_{T-\delta
}(\omega))-\varepsilon (C_{\delta }(\omega),R_{\delta })(\omega )
\end{equation*}%
where, by using (\ref{theta}),
\begin{align*}
\varepsilon (C_{\delta },R_{\delta })(\omega)
&\leq\frac{c_{d}}{\sqrt{\det C_{\delta }}}
(1+\Vert C_{\delta }^{-1/2}R_{\delta }\Vert
_{\delta ,2,q_{d}})^{\ell_{q}}\Vert C_{\delta }^{-1/2}R_{\delta
}\Vert _{\delta ,2,q_{d}} \\
&= \frac{c_{d}}{\sqrt{\det C_{\delta }}}(1+\theta _{\delta
,q_{d}})^{\ell_{q}}\theta _{\delta ,q_{d}}.
\end{align*}%
If $\omega \in \Gamma _{\delta ,r}(z)$ then $\theta _{\delta ,q_{d}}\leq
a_de^{-r^2}\leq 1$ so that%
$$
\varepsilon (C_{\delta },R_{\delta })(\omega )\leq \frac 12\times\frac{1}{(2\pi)^{d/2}\sqrt{\det
C_{\delta }}}\,e^{-r^2}.
$$
For $\omega \in \Gamma _{\delta ,r}(z)$ we also have
\begin{equation*}
\left\langle C_{\delta }^{-1}(F_{T-\delta }-z),F_{T-\delta }-z\right\rangle
=\left\vert F_{T-\delta }-z\right\vert _{\delta }^{2}\leq r^{2}
\end{equation*}%
so that%
\begin{equation*}
g_{C_{\delta }}(z-F_{T-\delta })\geq \frac{1}{(2\pi )^{d/2}\sqrt{\det
C_{\delta }}}e^{-r^{2}}.
\end{equation*}%
Then, by the choice of $a_de^{-r^2}$ we obtain
$$
\overline{p}_{F,\delta }(\omega ,z)\geq \frac{1}{(2\pi )^{d/2}\sqrt{\det
C_{\delta }}}e^{-r^{2}}-\varepsilon (C_{\delta },R_{\delta })(\omega )\geq
\frac{1}{2(2\pi )^{d/2}\sqrt{\det C_{\delta }}}e^{-r^{2}}.
$$
We conclude that
$$
\E(f(F)\mid \mathcal{F}_{T-\delta })(\omega )\geq \frac{1}{2e^{r^{2}}(2\pi
)^{d/2}}\int f(z)(\det C_{\delta })^{-1/2}1_{\Gamma _{\delta ,r}(z)}dz.
$$
$\Box$

We are now ready for our main result. It involves the concept of ``local densities'', that we define as follows: we say that the law of a r.v. $F$ taking values on $\R^d$ admits a local density around $y\in \R^d$ if there exists an open neighborhood $V_y$ of $y$ such that the restriction of the law of $F$ in $V_y$ is absolutely continuous w.r.t. the Lebesgue measure $\leb_d$ on $\R^d$. So, we have:

\begin{theorem}\label{th-pos}
For $\delta\in (0,T]$, let decomposition (\ref{deco}) hold and for $y\in\R^d$, $r>0$, assume that
$$
\P(\widetilde\Gamma_{\delta,r}(y))>0,
$$
where $\widetilde\Gamma_{\delta,r}(y)$ is the set in (\ref{gamma}). Then there exists $\eta>0$ and $c(y)>0$ such that for every Borel measurable set $A\subset B_\eta(y)$ one has
$$
\P(F\in A)\geq c(y)\leb_d(A).
$$
As a consequence, if the law of $F$ admits a local density $p_F$  around $y$ then one has
$$
p_F(x)\geq c(y)\quad \mbox{for a.e. $x\in B_\eta(y)$.}
$$
\end{theorem}

\def\ep{\varepsilon}
\textbf{Proof}.
For $\ep>0$, set
$$
\widetilde\Gamma_{\delta,r,\ep}(y)=\big\{\vert F_{T-\delta }-y\vert _{\delta }\leq
r/2\big\}\cap \big\{\det C_\delta\geq \ep\big\}\cap \big\{\theta _{\delta ,q_d}\leq
a_d\,e^{-r^2}\big\}.
$$
If $\P(\widetilde\Gamma_{\delta,r}(y))>0$ then there exists $\ep>0$ such that $\P(\widetilde\Gamma_{\delta,r,\ep}(y))>0$. On the set $\{\det C_\delta\geq \ep\}$, one has
$$
|\xi|_\delta\leq \ep^{-d/2}|\xi|,\quad\xi\in\R^d.
$$
Taking $\eta=\ep^{d/2}r/2$, one immediately has
$$
\widetilde\Gamma_{\delta,r,\ep}(y)\subset
\Gamma_{\delta,r}(x)\quad\mbox{for every $x\in B_\eta(y)$}
$$
where $\Gamma_{\delta,r}(x)$ is the set in (\ref{gamma0}). Therefore, by applying Lemma \ref{lemma-pos}, for every measurable and bounded function $f$ whose support is included in $B_\eta(y)$ one has
$$
\E(f(F)\,|\,\cl F_{T-\delta})(\omega)
\geq \frac{1}{2e^{r^2}(2\pi)^{d/2}}\,(\det C_\delta)^{-1/2}\, \I_{\tilde\Gamma_{\delta,r,\ep}(y)}\int_{B_\eta(y)} f(x)dx,
$$
and by passing to the expectation one gets the result with
$$
c(y)=\frac{1}{2e^{r^2}(2\pi)^{d/2}}\,\E((\det C_\delta)^{-1/2}\I_{\tilde\Gamma_{\delta,r,\ep}(y)})>0.
$$
$\Box$

\section{Examples}\label{sect-examples}
We apply now Theorem \ref{th-pos} to two cases in which a support theorem is available
and we give results for the strict positivity and lower bounds for the density which involve
suitable local or global non degeneracy conditions on the skeleton.

\subsection{It\^o processes}\label{sect-Ito}

We consider here a process $Z_t=(X_t,Y_t)^*$, taking values on $\R^d\times\R^n$,
which solves the following stochastic differential equation:
as $t\leq T$,
\begin{equation}\label{Ito}
\begin{array}{rl}
X_t=
&
\displaystyle
x_0+\sum_{j=1}^m\int_0^t\sigma_j(X_t,Y_t)dW^j_t+\int_0^tb(X_t,Y_t)dt\smallskip\\
Y_t=
&
\displaystyle
y_0+\sum_{j=1}^m\int_0^t\alpha_j(X_t,Y_t)dW^j_t+\int_0^t\beta(X_t,Y_t)dt.
\end{array}
\end{equation}
We are interested in dealing with strict positivity and/or lower bounds
for the probability density function of one component at a fixed time, say $X_T$,
as a consequence of  Theorem \ref{th-pos}.
This is a case in which a support theorem is available,
and we are going to strongly use it. For diffusion processes (that is, if we deal with $Z_T$ and not with $X_T$ only), we get an example which
is essentially the same as in the
paper of Ben Arous and L\'eandre \cite{bib:[BA.L]} and in the paper of Aida, Kusouka
and Stroock \cite{bib:[A.K.S]}.
Concerning the lower bounds, we will use lower estimates
for the probability that It\^o processes remains in a tube around a path proved in
Bally, Fern\'andez and Meda in \cite{bfm}.

\smallskip

So, in (\ref{Ito}) we assume that $\sigma_j,b\in C^4_b(\R^{d+n};\R^d)$
and $\alpha_j,\beta\in C^4_b(\R^{d+n};\R^n)$, $j=1,\ldots,m$, which implies that
$X^\ell_t,Y^i_t\in \D^{2,\infty}$ for all $\ell$ and $i$.

For $\phi\in L^2([0,T];\R^m)$, let $z_t(\phi)=(x_t(\phi),y_t(\phi)^*$
denote the skeleton associated to (\ref{Ito}), i.e.
\begin{equation}\label{skeleton}
\begin{array}{rl}
x_t(\phi)=
&
\displaystyle
x_0+\sum_{j=1}^m\int_0^t\sigma_j\big(x_t(\phi),y_t(\phi)\big)\phi^j_t\,dt
+\int_0^t\overline{b}\big(x_t(\phi),y_t(\phi)\big)dt\smallskip\\
y_t(\phi)=
&
\displaystyle
y_0+\sum_{j=1}^m\int_0^t\alpha_j\big(x_t(\phi),y_t(\phi)\big)\phi^j_t\,dt
+\int_0^t\overline{\beta}\big(x_t(\phi),y_t(\phi)\big)dt,
\end{array}
\end{equation}
in which
$\overline{b}=b-\frac 12 \sum_{j=1}^m\partial_{\sigma_j}\sigma_j$ and
$\overline{\beta}=\beta-\frac 12 \sum_{j=1}^m\partial_{\alpha_j}\alpha_j$,
where we have used the notation $(\partial_gf)^i=\<\nabla f^i, g\>$.

\smallskip

For a fixed $x\in\R^d$, we set
\begin{equation}\label{Cx}
\cl C(x)=\{\phi\in L^2([0,T];\R^m)\,:\, x_T(\phi)=x\}.
\end{equation}
We finally consider the following set of functions: for fixed $\mu\geq 1$ and $h>0$,
\begin{equation}\label{Lmuh}
L(\mu,h)=\{f\,:\,[0,T]\to\R_+\,;\, f_t\leq \mu f_s\mbox{ for all $t,s$ such that  $|t-s|\leq h$}\}.
\end{equation}

We have

\begin{theorem}\label{th-Ito}
Let $Z=(X,Y)^*$ denote the solution to (\ref{Ito}), with
$\sigma_j,b\in C^4_b(\R^{d+n};$ $\R^d)$
and $\alpha_j,\beta\in C^4_b(\R^{d+n};\R^n)$, $j=1,\ldots,m$.
Let $x\in \R^d$ be fixed and suppose that $\cl C(x)\neq \emptyset$. For $\phi\in \cl C(x)$, let $z_t(\phi)=(x_t(\phi),y_t(\phi))^*$ be as in (\ref{skeleton}).

\medskip

$\ i)$ Suppose there exists $\phi\in{\cal{C}}(x)$ such that $\sigma\sigma^*(x,y_T(\phi))>0$. Then there exists $\eta>0$ and $c(x)>0$ such that for every Borel measurable set $A\subset B_\eta(x)$ one has
$$
\P(X_T\in A)\geq c(x)\leb_d(A).
$$
In particular, if $X_T$ admits a local density $p_{X_T}$ around $x$ then $p_{X_T}\geq c(x)>0$ a.e. on the ball $B_\eta(x)$.

\smallskip

$\ ii)$
Suppose there exists $\phi\in\cl C(x)$ such that
$|\partial x_t(\phi)|\in L(\mu, h)$, for some $\mu\geq 1$ and $h>0$, and
$$
\sigma\sigma^* \big(x_t(\phi),y\big)\geq \lambda_*>0\quad
\mbox{for all $t\in [0,T]$ and $y\in \R^n$}.
$$
Then if the law of $X_T$ admits a continuous local density
$p_{X_T}$ around $x$ one has
$$
p_{X_T}(x)\geq
\Upsilon \exp\Big[-Q\Big(\Psi+
\frac{1}{\lambda_*}\int_0^{T}
|\partial_t x_t(\phi)|dt\Big)\Big],
$$
where $\Upsilon$, $Q$, $\Psi$
are all positive constants depending on
$d,T,\mu,h,\lambda_*$ and vector fields $\sigma_j,\alpha_j$, $j=1,\ldots,m$, and $b,\beta$.
\end{theorem}

In next Proposition \ref{prop-Ito} we study the existence of a local density and we prove in particular that under the requirement in part $ii)$, the local density really exists.
Actually, a little bit more work would show that the non degeneracy condition (\ref{Law1}) holds and by Lemma \ref{1} the local density is indeed continuous. But we are not interested here to enter in these technical arguments.

\medskip

\textbf{Proof of Theorem \ref{th-Ito}.}
$i)$ We take $0<\delta\leq T$ and we consider the decomposition $X_T=X_{T-\delta}+G_\delta+R_\delta$, where
\begin{align*}
G_\delta&=\sum_{j=1}^m\int_{T-\delta}^T\sigma_j(X_{T-\delta}, Y_{T-\delta})dW^j_t\\
R_\delta&=\sum_{j=1}^m\int_{T-\delta}^T\big(\sigma_j(X_t,Y_t)-\sigma_j(X_{T-\delta}, Y_{T-\delta})\big)dW^j_t
+\int_{T-\delta}^Tb(X_t,Y_t)dt.
\end{align*}
Conditionally on $\cl F_{T-\delta}$, the covariance matrix of the Gaussian r.v. $G_\delta$ is
$$
C_\delta=\sigma\sigma^*(X_{T-\delta}, Y_{T-\delta})\delta.
$$
So, we are in the framework studied in Section \ref{sect-perturbation}
and we proceed in order to apply Theorem \ref{th-pos}: $i)$ is proved as soon as  we find $\delta, r>0$ such that $\P(\widetilde\Gamma_{\delta,r}(x))>0$.

\smallskip

For $\phi\in \cl C(x)$, we denote $z^\phi(x)=(x,y_T(\phi))$ and we take $\phi$
such that $\sigma\sigma^*(z^\phi(x))>0$. We denote by $\lambda_*>0$ the lower eigenvalue of $\sigma\sigma^*(z^\phi(x))$.
Then, there exists $\ep>0$ such that
$$
\sigma\sigma^*(z)\geq \frac{\lambda_*}2\,I_d\quad\mbox{for every $z$ such that $|z-z^\phi(x)|<\ep$}.
$$
For a fixed $\delta\in(0,T]$, we have
$|z^\phi(x)-z_{T-\delta}(\phi)|
=|z_{T}(\phi)-z_{T-\delta}(\phi)|
\leq C(1+\|\phi\|_2)\sqrt\delta
=C_\phi\sqrt\delta$, so that
if $|Z_{T-\delta}-z_{T-\delta}(\phi)|<C_\phi\sqrt\delta$
then $|Z_{T-\delta}-z^\phi(x)|<2C_\phi\sqrt\delta$. We choose $\delta_0$ such that
$2C_\phi\sqrt\delta<\ep$ for all $\delta<\delta_0$.
So, if  $|Z_{T-\delta}-z_{T-\delta}(\phi)|<C_\phi\sqrt\delta$
we get
$$
C_\delta\geq \frac{\lambda_*}2\,\delta I_d
$$
and in particular,
\begin{align*}
|X_{T-\delta}-x|_\delta
&=|C_\delta^{-1/2}(X_{T-\delta}-x)|
\leq \Big(\frac2{\lambda_*\delta}\Big)^{1/2}\,
|X_{T-\delta}-x|\\
&\leq \Big(\frac2{\lambda_*\delta}\Big)^{1/2}\,
|Z_{T-\delta}-z^\phi(x)|
< \frac{2\sqrt 2\,C_\phi}{\sqrt{\lambda_*}}=:\frac r2.
\end{align*}
Moreover, for $q\geq 2$, a standard reasoning gives%
\begin{align*}
\|R_\delta\|^q_{\delta,2,q}
&=
\E(\left\vert R_{\delta }\right\vert ^{q}\mid
\mathcal{F}_{T-\delta })
+\E\Big(\sum_{l=1}^{2}\Big(\int_{[T-\delta
,T]^{l}}\left\vert D_{s_{1}...s_{l}}^{(l)}R_{\delta }\right\vert
^{2}ds_{1}...ds_{l}\Big)^{q/2}\mid \mathcal{F}_{T-\delta }\Big)\cr
&\leq
(C_{1,q}\delta)^q,
\end{align*}
so that
\begin{equation}\label{theta-th-Ito}
\theta_{\delta,q}
=\|C_\delta^{-1/2}R_\delta\|_{\delta,2,q}
\leq \frac 1{\sqrt{\lambda_*\delta}}\,
\|R_\delta\|_{\delta,2,q}
\leq
C_{2,q}\sqrt\delta.
\end{equation}
We take $\delta<\delta_0$ in order that $C_{2,q}\sqrt\delta<a_de^{-r^2}$.
For such a $\delta$ we get that $ \{|Z_{T-\delta}-z_{T-\delta}(\phi)|<C_\phi\sqrt\delta\}\subset \widetilde\Gamma_{\delta,r}(x)$ and by the classical support theorem for diffusion processes (see e.g. \cite{bib:IW}) one has $\P(\vert
Z_{T-\delta }-z_{T-\delta }(\phi )\vert < C_\phi\sqrt\delta)>0$, so that $\P(\widetilde \Gamma_{\delta,r}(x))>0$.

\smallskip

$ii)$
For $\xi\,:\,[0,T]\to\R^d$ and $R>0$, we set
$$
\tau^\phi_R(\xi)=\inf\{t\,:\,|\xi_t-x_t(\phi)|\geq R\}.
$$
We know that there exists $\phi\in \cl C(x)$ and
$\ep>0$ such that if $\tau^\phi_\ep(\xi)>T$ then
$$
\sigma\sigma^* (\xi_t,y)\geq \lambda_*I_d
$$
for any $t\in [0,T]$ and $y\in \R^n$. So,
on the set $\{\tau^\phi_\ep(X)>T\}$ one gets $C_\delta\geq $ $\lambda_*\delta I_d$. Moreover,
if $\tau^\phi_\ep(X)>T$ then for $0<\delta<T$
\begin{align*}
|X_{T-\delta}-x|
&=|X_{T-\delta}-x_T(\phi)|
\leq |X_{T-\delta}-x_{T-\delta}(\phi)|+|x_{T-\delta}(\phi)-x_T(\phi)|\\
&< \ep+\int_{T-\delta}^T|\partial_t x_t(\phi)|dt.
\end{align*}
Since again (\ref{theta-th-Ito}) holds,
we take $\delta<T $ such that $\int_{T-\delta}^T|\partial_t x_t(\phi)|dt<\ep$
and $\theta_{\delta,q_d}\leq $ $a_de^{-(2\ep)^2}$. Therefore,
$\{\tau^\phi_\ep(X)>T\}\subset \Gamma_{\delta,2\ep}(x)$
and by using Lemma \ref{lemma-pos} we get
$$
p_{X_T}(x)\geq \frac 1{2(2\pi\lambda_*\delta)^{d/2}e^{4\ep^2}}\,
\P\big(\tau^\phi_\ep(X)>T\big)\equiv
 \Upsilon\times \P\big(\tau^\phi_\ep(X)>T\big) .
$$
Now, the hypotheses allow one to use Theorem
1 in Bally,
Fern\'andez and Meda \cite{bfm}: one has
$$
\P(\tau^\phi_{\ep}(X)>T)\geq  \exp\Big(-Q\Big(\Psi+
\frac{1}{\lambda_*}\int_0^{T}
|\partial_t x_t(\phi)|dt\Big)\Big)
$$
and the statement holds.
$\square$

\begin{example}\label{grushin}
Let $n\geq 1$ and $k\geq 0$ be fixed integers and let $(X,Y)$ be the $2$-dimensional process solution to
\begin{align*}
&X_t=x_0+\int_0^t Y_s^n\,dW^1_s+\int_0^tY_s^k\,ds,\\
&Y_t=y_0+W^2_t,
\end{align*}
$W$ denoting a Brownian motion in $\R^2$. The pair $(X,Y)$ then follows the well-known Grushin diffusion. Here, we are interested  in the study of the component $X$ only, because this gives an example in between the two cases studied in Theorem \ref{th-Ito}. In fact, one has $\sigma\sigma^*(x,y)=y^{2n}$, and this vanishes as $y=0$, so there is no hope that part $ii)$ holds. Nevertheless, $i)$ is always true. In fact, since the strong H\"ormander condition holds for the diffusion pair $(X,Y)$, the law of $(X_T,Y_T)$ has a smooth density in $\R^2$, so that $X_T$ has a smooth density as well. Moreover, the associated skeleton is given by
$$
\begin{array}{rl}
x_t(\phi)=
&
\displaystyle
x_0+\int_0^ty_t^n(\phi)\phi^1_t\,dt
+\frac 12\,\int_0^t\big(2y_t^k(\phi)-n\,y_t^{2n-1}(\phi)\big)dt\smallskip\\
y_t(\phi)=
&
\displaystyle
y_0+\int_0^t\phi^2_t\,dt,
\end{array}
$$
so it is clear that for every $x\in\R$ one has ${\cal{C}}(x)\neq\emptyset$ and one can choose $\phi\in{\cal{C}}(x)$ such that $\sigma\sigma^*(x,y_T(\phi))>0$, that gives $p_{X_T}(x)>0$.

\end{example}

We propose now a sufficient condition for the existence of a local density, that in particular says that under the hypothesis of $ii)$ in Theorem \ref{th-Ito}, a local density really exists.

\begin{proposition}\label{prop-Ito}
Set
$$
{\cal O} =\{x\in\R^d\,:\,\P(\sigma \sigma ^{\ast }(x,Y_{T})>0)=1\}.
$$
Then for every $x\in {\cal O}$ the law of $X_T$ admits a local density $p_{X_T}$ around $x$. As a consequence, if $x\in {\cal O}$ is such that ${\cal C}(x)\neq \emptyset$ and for some $\phi\in{\cal{C}}(x)$ one has $\sigma\sigma^*(x,y_T(\phi))>0$, then the local density $p_{X_T}$ is a.e. strictly positive around $x$.
\end{proposition}

\textbf{Proof.} For $x\in {\cal O}$, set $D_{x}=\{y\in\R^n\,:\,\sigma \sigma ^{\ast
}(x,y)>0\}$. $D_x$ is an open set, so there exist a sequence $\{y_{i}\}\subset \R^n$ and a sequence $\{r_{i}\}_i\subset \R_+$ such that
$D_{x}=\cup _{i\in \N}B_{\frac{1}{2}r_{i}}(y_{i})$ and $B_{r_{i}}(y_{i})%
\subset D_{x}.$ Moreover $\sigma \sigma ^{\ast }(\bar x,y)\geq \lambda
_{i}>0$ for every $y\in B_{r_{i}}(y_{i})$ and $\bar x\in B_{r}(x).$
For any fixed $i$, we consider a localizing r.v. $U_i$ of the form (\ref{Mall10'}): we set
$$
U_i=\psi_{r^2} (\left\vert X_{T}-x\right\vert^2 )\psi _{r_i^2}(\left\vert
Y_{T}-y_{i}\right\vert^2 ).
$$
By (\ref{Mall11}), $U_i$ is a good localizing r.v. (that is, (\ref{Mall1}) holds), and we set $d\P_i=U_id\P$. We also notice that if $U_i\neq 0$ then $\sigma \sigma ^{\ast
}(X_{T},Y_{T})\geq \lambda _{i}>0$. This property allows one to use a standard argument showing that the Malliavin covariance matrix of $F=X_{T}$ has finite inverse moments of any
order with respect to $\P_{i}$, which means that (\ref{Law1}) holds. So, we can use Lemma
\ref{1} and we can conclude that the law of $X_{T}$ with
respect to $\P_{i}$ is absolutely continuous with respect to the Lebesgue
measure. Take now $A\subset B_{r/2}(x)$ a set of Lebesgue measure equal to zero.
Since $x\in {\cal O} $ we have $\P(Y_T\in D_{x})=1$ so
\begin{align*}
\P(X_{T} \in A)
&=\P(X_{T}\in A,Y_{T}\in D_{x})\leq \sum_{i}\P(X_{T}\in
A,Y_{T}\in B_{\frac{1}{2}r_{i}}(y_{i})) \\
&=\sum_{i}\P_{i}(X_{T}\in A,Y_{T}\in B_{\frac{1}{2}r_{i}}(y_{i}))
\end{align*}%
the last equality being true because
$\psi_{r^2} (\left\vert X_{T}-x\right\vert^2 )\psi _{r_i^2}(\left\vert
Y_{T}-y_{i}\right\vert^2 )=1$ if $X_{T}\in A$ and $Y_{T}\in B_{\frac{1}{2}r_{i}}(y_{i})$.
Since the law of $X_{T}$ under $\P_{i}$ is absolutely continuous
with respect to the Lebesgue measure we obtain $\P_{i}(X_{T}\in A,Y_{T}\in
B_{\frac{1}{2}r_{i}}(y_{i}))=0$ for every $i$, and this proves that a local density $p_{X_T}$ around $x$ exists. The final statement comes now immediately from Theorem \ref{th-Ito}.
$\Box$

\begin{example}\label{new-ex}
Consider the diffusion process
\begin{align*}
X_{t}^{1} &=x^{1}_0+\int_{0}^{t}\alpha (\left\vert X_{s}\right\vert
)\left\vert Y_{s}\right\vert \circ dW_{s}^{1}+\int_{0}^{t}\left\vert
X_{s}\right\vert \circ dW_{s}^{3}, \\
X_{t}^{2} &=x^{2}_0+\int_{0}^{t}\alpha (\left\vert X_{s}\right\vert
)\left\vert Y_{s}\right\vert \circ dW_{s}^{2}+\int_{0}^{t}\left\vert
X_{s}\right\vert \circ dW_{s}^{3}, \\
dY_{t} &=y_0+\int_{0}^{t}\beta (X_{s})\circ dW_{t}^{4}
\end{align*}%
where $W$ is a standard Brownian motion taking values in $\R^4$ and $\alpha ,\beta $
are $C_{b}^{4}$ functions. We suppose that $\{r:\alpha (r)\neq
0\}=B_1(0)$ and that $\beta(x_0)\neq 0$, the latter requirement ensuring in particular that the law of $Y_T$ has a density. Therefore, for every $x\in B_1(0)$ one has $\P(\sigma\sigma^*(x,Y_T)>0)=1$ and by applying  Proposition \ref{prop-Ito} one gets that $X_T$ has a local density around every point in $B_1(0)$.
Now, in order to study its positivity property, let us write down the associated skeleton:
for a square integrable  control path $\phi$, one has
\begin{align*}
x_{t}^{1}(\phi ) =&x^{1}_0+\int_{0}^{t}\alpha (\left\vert x_{s}(\phi
)\right\vert )\left\vert y_{t}(\phi )\right\vert \phi
_{s}^{1}ds+\int_{0}^{t}\left\vert x_{s}(\phi )\right\vert \phi _{s}^{3}ds \\
x_{t}^{2}(\phi ) =&x^{2}_0+\int_{0}^{t}\alpha (\left\vert x_{s}(\phi
)\right\vert )\left\vert y_{t}(\phi )\right\vert \phi
_{s}^{2}ds+\int_{0}^{t}\left\vert x_{s}(\phi )\right\vert \phi _{s}^{3}ds, \\
y_{t}(\phi ) =&y_0+\int_{0}^{t}\beta (x_{s}(\phi )))\phi _{s}^{4}ds.
\end{align*}%
We recall that the support theorem of Stroock and Varadhan
asserts that the law of $(X_{t}^{1},X_{t}^{2},Y_{t})_{t\geq 0}$ is the
closure (with respect to the uniform norm) of the points of the skeleton as above.
Notice that if $\left\vert x_{t}^{1}(\phi )\right\vert \geq 1$ then $\alpha (\left\vert
x_{t}(\phi )\right\vert )=0$ and so $\partial _{t}x_{t}^{1}(\phi
)=\left\vert x_{t}(\phi )\right\vert \phi _{t}^{3}=\partial
_{t}x_{t}^{2}(\phi ).$ This means that outside the unit ball the skeleton $%
(x_{t}^{1}(\phi ),x_{t}^{2}(\phi ))$ may travel on a line which is parallel to the principal diagonal (i.e. $x^1=x^2$), but only on this line. If $\left\vert x_{t}^{1}(\phi )\right\vert <1$ then one may use the controls $\phi ^{1}$ and $\phi ^{2}$ and then $(x_{t}^{1}(\phi ),x_{t}^{2}(\phi ))$ may travel in any direction inside the open unit ball. Having this in mind, we
define the strip $S=\{(x^1,x^2)\,:\,|x^1-x^2|<\sqrt 2
\}$ and thanks to the above discussion and the support theorem we have the following three
cases.

\begin{description}
\item[{\rm 1.}] {\rm $x_0\notin S$.}
Here, for every $s$ the law of $X_{s}$ is concentrated on
the line which is parallel to the principal diagonal and contains $x_0$. In particular, $\alpha(|X_s|)=0$ a.s. for every $s$, so $X$ is actually a diffusion process  satisfying
\begin{equation*}
X_{t}^{1}=x^{1}_0+\int_{0}^{t}\left\vert X_{s}\right\vert \circ dW_{s}^{3},\quad
X_{t}^{2}=x^{2}_0+\int_{0}^{t}\left\vert X_{s}\right\vert \circ  dW_{s}^{3}.
\end{equation*}%

\item[{\rm 2.}] {\rm $x_0\in S$ but $x_0\notin B_1(0)$.}
Here the support of the law of $X_{T}$ is the whole $S$. By using Proposition \ref{prop-Ito}, we can say that $X_T$ has a local density around any point in  $B_1(0)$ and moreover, there exists a version of the local density which is strictly positive in the ball. But we have no information outside the ball.

\item[{\rm 3.}] {\rm $x_0\in B_1(0)$.}
We can assert the same statements as in case 2 but with some refinements. In fact, here if $y_0\neq 0$ then $\alpha(x_0)y_0\neq 0$, so that the law of $X_T$ has a smooth global density which is strictly positive on the unit ball $B_1(0)$.
\end{description}

Concerning point $ii)$ of Theorem \ref{th-Ito}, it does not apply except when $x_0,x\in B_1(0)$.

\end{example}

\subsection{Diffusion processes satisfying a weak H\"{o}rmander
condition: an example}\label{asian}
In this section we treat an example of diffusion process which satisfies the
weak H\"{o}rmander condition and has been recently studied in Bally and Kohatsu-Higa
\cite{bib:[B.KH]} (we are going to use the ideas and the
estimates from that paper).
Since lower bounds for the density have been already discussed in \cite{bib:[B.KH]},
we deal here only with the strict positivity. So, we give an application of our Theorem
\ref{th-pos} in a case of degenerate diffusion coefficients.

\smallskip

We consider the diffusion process
\begin{equation}
X_{t}^{1}=x^{1}+\int_{0}^{t}\sigma_1
(X_{s})dW_{s}+\int_{0}^{t}b_{1}(X_{s})ds,\quad
X_{t}^{2}=x^{2}+\int_{0}^{t}b_{2}(X_{s})ds  \label{H1}
\end{equation}%
and we assume that $\sigma_1 ,b_{1},b_{2}\in C_{b}^{\infty }(\R^{2};\R).$ Actually,
it suffices that they are four times differentiable - but we do not focus on
this aspect here. Moreover, we fix some point $y\in \R^{2}$ and we assume
that
\begin{equation}
\left\vert \sigma_1 (y)\right\vert > c_{\ast }>0\quad\mbox{and}\quad \left\vert \partial
_{1}b_{2}(y)\right\vert > c_{\ast }>0.  \label{H2}
\end{equation}%
Let $\sigma =(\sigma_1 ,0)^*$ and $b=(b_{1},b_{2})^*.$ The Lie bracket $%
[\sigma ,b]$ is computed as%
$$
[\sigma,b](x) =\partial _{\sigma%
}b(x)-\partial _{b}\sigma(x)=\binom{\sigma_1 (x)\partial _{1}b_{1}(x)-b_{1}(x)\partial _{1}\sigma_1
(x)-b_{2}(x)\partial _{2}\sigma_1 (x)}{\sigma_1 (x)\partial _{1}b_{2}(x)}.
$$
So assumption (\ref{H2}) is equivalent with the fact that $\sigma%
(y)$ and $[\sigma,b](y)$ span $\R^{2},$ and this is the weak H%
\"{o}rmander condition in $y.$

We set $\overline{b}=b-\frac 12\partial_\sigma\sigma $ and
for a measurable function $\phi\in L^2([0,T],\R)$
we consider the skeleton $x(\phi )$, i.e. the solution to the equation%
\begin{equation*}
x_{t}(\phi )=x+\int_{0}^{t}\Big(\sigma (x_{s}(\phi ))\phi _{s}+\overline{b%
}(x_{s}(\phi ))\Big)ds.
\end{equation*}

\begin{proposition}\label{hor-prop}
Assume that $\sigma_1 ,b_{1},b_{2}\in C_{b}^{\infty }(\R^{2})$ and (\ref{H2})
holds. Then the law of $X_{T}$ has a local smooth density $p_{T}(x,\cdot )$\ in a
neighborhood of $y.$ Moreover, if there exists
a control $\phi\in L^2([0,T]) $
such that $x_{T}(\phi )=y$ then $p_{T}(x,y)>0.$
\end{proposition}

Before starting with the proof of Proposition \ref{hor-prop}, let us
consider the following
decomposition: for $\delta\in (0,T]$, we set
\begin{equation}\label{hor-dec}
F=X_T-x_T(\phi)\quad\mbox{ and } \quad
F=F_{T-\delta}+G_\delta+R_\delta
\end{equation}
where $F_{T-\delta}=X_{T-\delta}-x_{T-\delta}(\phi)$ and
\begin{align*}
G_\delta^1
&=\int_{T-\delta}^T \sigma_1(X_{T-\delta})dW_s,
\quad G_\delta^2
=\int_{T-\delta}^T \partial_\sigma b_2(X_{T-\delta})(T-s)dW_s\\
R^1_\delta
&=\int_{T-\delta}^T\Big(\sigma_1(X_s)-\sigma_1(X_{T-\delta})\Big)dW_s
+\int_{T-\delta}^Tb_1(X_s)ds+\\
&\quad
-\int_{T-\delta}^T\Big(\sigma_1\big(x_s(\phi)\big)\phi_s+\overline{b}_1\big(x_{T-\delta}(\phi)\big)\Big)ds\\
R^2_\delta
&=\int_{T-\delta}^T\Big(\partial_\sigma b_2(X_s)-\partial_\sigma b_2(X_{T-\delta})\Big)(T-s)dW_s
+\delta\Big(b_2(X_{T-\delta})-b_2\big(x_{T-\delta}(\phi)\big)\Big)+\\
&\quad +\int_{T-\delta}^TLb_2(X_s)(T-s)ds
-\int_{T-\delta}^T\Big(b_2\big(x_s(\phi)\big)-
b_2\big(x_{T-\delta}(\phi)\big)\Big)ds
\end{align*}
in which
$L=\frac{1}{2}\sigma \sigma ^{\ast }\partial _{x}^{2}+b\partial _{x}$ denotes the infinitesimal generator of $X$.

The covariance matrix of the conditional (on $\cl F_{T-\delta}$) Gaussian r.v. $G_\delta$
is given by
$$
C_\delta=\delta\sigma_1^2(X_{T-\delta})\left(
\begin{array}{cc}
1 &
\partial_1b_2(X_{T-\delta})\frac {\delta}2\smallskip\\
\partial_1b_2(X_{T-\delta})\frac {\delta}2 &
(\partial_1b_2)^2(X_{T-\delta})\frac {\delta^2}3
\end{array}
\right).
$$
We need now some estimates which can be easily deduced from \cite{bib:[B.KH]}. In order
to be self contained, we propose here the following
\begin{lemma}\label{hor-lemma}
Let $\rho^2_\delta=\max(\delta,\int_{T-\delta}^T|\phi_s|^2ds)$.
Then, there exist $\delta_0>0$ such that for every $\delta<\delta_0$, on the set
$\{|F_{T-\delta}|<\delta^{3/2}\rho_\delta\}$ the following properties hold:

\begin{itemize}
\item[$i)$]
$\displaystyle\det C_\delta\geq c_1\frac{\delta^4}{12}$;
\item[$ii)$]
for every $\xi\in \R^2$,
$\displaystyle
|\xi|^2_\delta\leq \frac{c_2}{\delta^3}\Big(\delta^2|\xi_1|^2+|\xi_2|^2\Big);
$
in particular,
$|F_{T-\delta}|_\delta\leq c_2\rho_\delta$;
\item[$ii)$]
 for every $q\geq 2$, $\theta_{\delta,q}\leq L_q \rho_\delta$.
\end{itemize}
Here, $c_1$, $c_2$  and $L_q$ are suitable positive constants
depending on $c_*$ and upper bounds for $\sigma$, $b$ and their derivatives up to order $4$,
$L_q$ depending on $q$ also, and we recall that $|\xi|_\delta=|C_\delta^{-1/2}\xi|$.

\end{lemma}

\textbf{Proof.}
First, by recalling that $\sigma$ and $b$ are bounded, for
some positive constant $C$ we have
$$
|x_{T-\delta}(\phi)-x_T(\phi)|
\leq C\rho^2_\delta
$$
so that
$$
|X_{T-\delta}-y|\leq |X_{T-\delta}-x_{T-\delta}(\phi)|+C\rho^2_\delta.
$$
Therefore, we can choose $\delta_0$ such that for all $\delta<\delta_0$
the following holds:
if $|F_{T-\delta}|=|X_{T-\delta}-x_{T-\delta}(\phi)|<\delta^{3/2}\rho_\delta$ then
$$
|\sigma_1(X_{T-\delta})|\geq c_*>0
\quad\mbox{and}\quad
|\partial_1b_2(X_{T-\delta})|\geq c_*>0.
$$
Therefore,
$$
\det C_\delta=\frac{(\sigma_1^2\partial_1b_2)^2(X_{T-\delta})\delta^4}{12}
\geq c_1\delta^4
$$
and $i)$ holds.
Moreover, we have
$$
C_\delta^{-1}=\frac{1}{(\sigma_1\partial_1b_2)^2(X_{T-\delta})\delta^3}
\left(
\begin{array}{cc}
4(\partial_1b_2)^2(X_{T-\delta})\delta^2 &
-6(\partial_1b_2)(X_{T-\delta})\delta\smallskip\\
-6(\partial_1b_2)(X_{T-\delta})\delta&
12
\end{array}
\right)
$$
so that for $\xi\in\R^2$,
\begin{align*}
|C_\delta^{-1/2}\xi|^2
&=\<C_\delta^{-1}\xi,\xi\>
=\frac{1}{(\sigma_1\partial_1b_2)^2(X_{T-\delta})\delta^3}\Big(
\big(2\partial_1 b_2(X_{T-\delta})\delta\xi_1-3\xi_2\big)^2+3\xi_2^2\Big)\\
&
\leq \frac{C}{c_*^4\delta^3}\Big(\delta^2|\xi_1|^2+
 |\xi_2|^2\Big)
\end{align*}
where $C$ depends on $\sigma$ and $b$. Then, if
$|F_{T-\delta}|<\delta^{3/2}\rho_\delta$ one gets
$$
|F_{T-\delta}|_\delta^2
=|C_\delta^{-1/2}F_{T-\delta}|^2
\leq \frac{C}{c_*^4\delta^3}\delta^3\rho_\delta^2(\delta^2+1)
\leq c_2 \rho_\delta^2
$$
and $ii)$ is proved.
As for $iii)$,  for $q\geq 2$ we have
$$
\E\big(|C_\delta^{-1/2}R_\delta|^q\big)
\leq \Lambda_q\Big(\E(|\delta^{-1/2}R^1_\delta|^q)+
\E( |\delta^{-3/2}R^2_\delta|^q)\Big),
$$
where $\Lambda_q$ depends on $q$, $c_*$, $\sigma$ and $b$.
Now, by using the Burkholder inequality
and the boundedness of the coefficients $b$ and $\sigma$ and of their
derivatives, one has
\begin{align*}
\E(|\delta^{-1/2}R^1_\delta|^q)
&\leq C_q\delta^{-q/2}\Big[
\E\Big(\Big|\int_{T-\delta}^T\big(\sigma_1(X_s)-\sigma_1(X_{T-\delta})\big)dW_s\Big|^q\Big)+
\\
&\qquad\qquad\qquad
+\E\Big(\Big|\int_{T-\delta}^Tb_1(X_s)ds\Big|^q\Big)+\\
&\qquad\qquad\qquad
+\E\Big(\Big|
\int_{T-\delta}^T\Big(\sigma_1\big(x_s(\phi)\big)\phi_s
+\overline{b}_1\big(x_{T-\delta}(\phi)\big)\Big)ds\Big|^q\Big)\Big]\\
&\leq C_qC\delta^{-q/2}\cdot\Big(\delta^{q}+\delta^{q/2}\Big(\int_{T-\delta}^T|\phi_s|^2ds\Big)^{q/2}\Big)
\leq 2C_qC\rho_\delta^{q}
\end{align*}
where $C_q$ depends on $q$ only and $C$ depends on the bounds of the diffusion coefficients.
Similarly (in the following $C$ denotes a suitable constant),
\begin{align*}
\E(|\delta^{-3/2}R^2_\delta|^q)
&\leq C_q\delta^{-3q/2}\Big[
\E\Big(\Big|\int_{T-\delta}^T\big(\partial_{\sigma} b_2(X_s)-\partial_{\sigma} b_2(X_{T-\delta})\big)(T-s)dW_s\Big|^q\Big)+\\
&\qquad\qquad\qquad
+\E\Big(
\delta^q\big|b_2(X_{T-\delta})-b_2\big(x_{T-\delta}(\phi)\big)\big|^q\Big)+\\
&\qquad\qquad\qquad
+\E\Big(\Big|\int_{T-\delta}^TL b_2(X_s)(T-s)ds\Big|^q\Big)+\\
&\qquad\qquad\qquad
+\Big|\int_{T-\delta}^T\Big(b_2\big(x_s(\phi)\big)-
b_2\big(x_{T-\delta}(\phi)\big)\Big)ds\Big|^q\Big]\\
&\leq 2C_qC\delta^{-3q/2}\Big(
\delta^{2q}+\delta^q|F_{T-\delta}|^q+\delta^q\sup_{T-\delta\leq s\leq T}|x_s(\phi)-x_{T-\delta}(\phi)|^q\Big)\\
&\leq 2C_qC\delta^{-3q/2}\Big(
\delta^{2q}+\delta^q\cdot\delta^{3q/2}\rho_\delta^q+
\delta^q\cdot\Big[\delta^q+\delta^{q/2}\cdot\Big(\int_{T-\delta}^T|\phi_s|^2ds\Big)^{q/2}\Big]\Big)\\
&\leq C_qC\rho_\delta^{q}.
\end{align*}
The same arguments may be used to give upper estimates for the remaining terms
in $\|C_\delta^{-1/2}R_\delta\|^q_{\delta,2,q}$ that contain
the Malliavin derivatives. So, we deduce that
$$
\|C_\delta^{-1/2}R_\delta\|_{\delta,2,q}\leq L_q\rho_\delta
$$
and the proof is completed. $\square$

\medskip

We are now ready for the

\smallskip

\textbf{Proof of Proposition \ref{hor-prop}.}
Consider the decomposition (\ref{hor-dec}): we have $p_{X_T}(y)= p_F(0)$.
We use Lemma \ref{hor-lemma} and Theorem \ref{th-pos}.
So, there exists $\delta_0$ such that for $\delta<\delta_0$
if $|F_{T-\delta}|<\delta^{3/2}\rho_\delta$ then $|F_{T-\delta}|_\delta<c_2\rho_\delta$.
We take now $\delta_1<\delta_0$ and $r=c_2\rho_{\delta_1}$. So, there exists $\delta<\delta_1$
such that $\theta_{\delta,q_d}=\|C_\delta^{-1/2}R_\delta\|_{\delta,2,q_d}\leq a_de^{-r^2}$. Therefore,
$\widetilde\Gamma_{\delta,r}(0)\supset \{|F_{T-\delta}|<\delta^{3/2}\rho_\delta\}$ and
by the support theorem one has
$\P(|F_{T-\delta}|<\delta^{3/2}\rho_\delta)=\P(|X_{T-\delta}-x_{T-\delta}(\phi)|<\delta^{3/2}\rho_\delta)>0$,
so Theorem \ref{th-pos} allows one to conclude.
$\square$

\appendix
\section{Proof of (\ref{1bisbis}) and (\ref{2bisbis})}\label{bisbis}
We recall inequalities (\ref{1bisbis}) and (\ref{2bisbis}): they are given by
\begin{equation}\label{1bis}
\left\vert \widehat\sigma _{F}^{i,j}-\widehat\sigma _{G}^{i,j}\right\vert
\leq
C(1\vee \det \widehat\sigma
_{F})(1\vee\det \widehat\sigma _{G})\left\vert D(F-G)\right\vert\big(1+\left\vert DF\right\vert +\left\vert
DG\right\vert \big)^{2d-1}
\end{equation}
and 
\begin{equation}\label{2bis}
\begin{array}{c}
\left\vert D\widehat\sigma _{F}^{i,j}-D\widehat\sigma _{G}^{i,j}\right\vert
\leq  C(1\vee \det \widehat\sigma
_{F})^2(1\vee\det \widehat\sigma _{G})^{2} (\left\vert D(F-G)\right\vert +\left\vert D^{(2)}(F-G)\right\vert )\times\smallskip\\
\times \big(1+\left\vert DF\right\vert +\left\vert
DG\right\vert +\left\vert D^{(2)}F\right\vert +\left\vert D^{(2)}G\right\vert
\big)^{6d-3}
\end{array}
\end{equation}
respectively.

\medskip

\textbf{Proof of (\ref{1bis}).}
First, one has
\begin{align*}
|\sigma_F^{ij}-\sigma_G^{ij}|
&= |\<DF^i,DF^j\>-\<DG^i,DG^j\>|\\
&\leq |DF^j|\,|DF^i-DG^i|+|DG^i|\,|DF^j-DG^j|\\
&\leq |DF-DG|(|DF|+|DG|),
\end{align*}
so that
\begin{equation}\label{a}
\sup_{i,j}|\sigma_F^{ij}-\sigma_G^{ij}|\leq |DF-DG|\big(|DF|+|DG|\big).
\end{equation}
Now, (\ref{a}) gives
\begin{equation}\label{b}
|\det \sigma_F-\det\sigma_G|\leq C\,|DF-DG|\big(|DF|+|DG|\big)^{2d-1}.
\end{equation}
In fact, setting $\cl P_d$ as the set of all the permutations of $(1,\ldots,d)$, by using (\ref{a}) one has
($C$ denotes a positive constant, depending on $d$, that may vary from line to line)
\begin{align*}
|\det\sigma_F-\det\sigma_G|
\leq &\sum_{\gamma\in\cl P_d}|\sigma_F^{1\gamma_1}\cdots\sigma_F^{d\gamma_d}
-\sigma_G^{1\gamma_1}\cdots\sigma_G^{d\gamma_d}|\\
\leq &
\sum_{\gamma\in\cl P_d}\Big(|\sigma_F^{1\gamma_1}-\sigma_G^{1\gamma_1}|\,
|\sigma_F^{2\gamma_2}\cdots\sigma_F^{d\gamma_d}|+|\sigma_G^{1\gamma_1}|\,
|\sigma_F^{2\gamma_2}\cdots\sigma_F^{d\gamma_d}
-\sigma_G^{2\gamma_2}\cdots\sigma_G^{d\gamma_d}|\Big)\\
\leq & C|DF-DG|(|DF|+|DG|)|DF|^{2(d-1)}+\\
&+|DG|^2\sum_{\gamma\in\cl P_d}|\sigma_F^{2\gamma_2}\cdots\sigma_F^{d\gamma_d}
-\sigma_G^{2\gamma_2}\cdots\sigma_G^{d\gamma_d}|\\
\leq & \cdots\\
\leq &
C|DF-DG|(|DF|+|DG|)\sum_{\ell=0}^{d-1}|DG|^{2\ell}|DF|^{2(d-1-\ell)}\\
\leq &
C|DF-DG|(|DF|+|DG|)\times\big(|DF|+|DG|\big)^{2(d-1)}
\end{align*}
and (\ref{b}) follows.

With $\sigma$ generally denoting $\sigma_F$ and/or $\sigma_G$, let us now set $\tilde \sigma$ as the matrix of cofactors, so that $\widehat\sigma^{ij}=(-1)^{i+j}(\det\sigma)^{-1}\tilde\sigma^{ji}$. By recalling that $\tilde\sigma^{ji}$ is the determinant of the sub-matrix of $\sigma$ obtaining by deleting the $j$th row and the $i$th column of $\sigma$, (\ref{b}) gives
\begin{equation}\label{c}
\sup_{i,j}|\tilde\sigma^{ij}_F-\tilde\sigma^{ij}_G|\leq C\,|DF-DG|\big(|DF|+|DG|\big)^{2d-3}.
\end{equation}
So,  we have
\begin{align*}
\vert \widehat\sigma _{F}^{i,j}-\widehat\sigma _{G}^{i,j}\vert
=&
\vert (\det\sigma_F)^{-1}\tilde\sigma _{F}^{ji}-(\det\sigma _{G})^{-1}\tilde\sigma_G^{ji}\vert\\
\leq &
(\det\sigma_F)^{-1}|\tilde\sigma_F^{ji}-\tilde\sigma_G^{ji}|
+(\det\sigma_F)^{-1}(\det\sigma_G)^{-1}|\tilde\sigma_G^{ji}|\,|\det\sigma_F-\det\sigma_G|\\
\leq &
(1\vee\det\widehat\sigma_F)(1\vee\det\widehat\sigma_F)
\Big(|\tilde\sigma_F^{ji}-\tilde\sigma_G^{ji}|+|\det\sigma_F-\det\sigma_G|\Big).
\end{align*}
Now, in the brackets of the r.h.s. we use (\ref{b}) and (\ref{c}), and we obtain
\begin{align*}
\Big\vert \widehat\sigma _{F}^{i,j}-\widehat\sigma _{G}^{i,j}\vert
\leq &
C\,(1\vee\det\widehat\sigma_F)(1\vee\det\widehat\sigma_G)
|DF-DG|\times\\
&\times\Big(\big(|DF|+|DG|\big)^{2d-3}+\big(|DF|+|DG|\big)^{2d-1}\Big),
\end{align*}
so that (\ref{1bis}) holds. $\square$

\medskip

\textbf{Proof of (\ref{2bis}).} The proof is similar to the one above. Here, it is sufficient to use the following further inequalities, straightforward to be proven:
\begin{equation}\label{d}
\begin{array}{rl}
\sup_{i,j}|D\sigma_F^{ij}-D\sigma_G^{ij}|
\leq &\big(|D(F-G)|+|D^{(2)}(F-G)|\big)\times\smallskip\\
&\times\big(|DF|+|DG|
+|D^{(2)}F|+|D^{(2)}G|\big)
\end{array}
\end{equation}
This implies that
\begin{equation}\label{f}
\begin{array}{rl}
|D\det \sigma_F-D\det\sigma_G|
\leq& C\,\big(|D(F-G)|+|D^{(2)}(F-G)|\big)\times\smallskip\\
&\times\big(|DF|+|DG|+|D^{(2)}F|+|D^{(2)}G|\big)^{2d-1}\smallskip\\
\end{array}
\end{equation}
and, as a consequence, 
\begin{equation}\label{g}
\begin{array}{rl}
\sup_{ij}|D\tilde\sigma^{ij}_F-D\tilde\sigma^{ij}_G|
\leq & C\,\big(|D(F-G)|+|D^{(2)}(F-G)|\big)\times\smallskip\\
&\times\big(|DF|+|DG|+|D^{(2)}F|+|D^{(2)}G|\big)^{2d-3}
\end{array}
\end{equation}
By developing $D\widehat\sigma _{F}^{i,j}=D\big((\det\sigma_F)^{-1}\tilde\sigma^{ji}_F\big)$ and $D\widehat\sigma _{G}^{i,j}=D\big((\det\sigma_G)^{-1}\tilde\sigma^{ji}_G\big)$ and by using also (\ref{d}), (\ref{f}) and (\ref{g}), one obtains (\ref{2bis}). $\square$

\end{document}